\documentclass[sigconf]{acmart}

\AtBeginDocument{%
  \providecommand\BibTeX{{%
    \normalfont B\kern-0.5em{\scshape i\kern-0.25em b}\kern-0.8em\TeX}}}

\usepackage{multirow}
\usepackage{xspace}
\newcommand{\method}{\texttt{MTC}\xspace }

\usepackage{amsmath,amsthm}
\theoremstyle{definition}

\newtheorem{problem}{Problem}

\usepackage{listings}
\definecolor{dkgreen}{rgb}{0,0.6,0}
\definecolor{gray}{rgb}{0.5,0.5,0.5}
\definecolor{mauve}{rgb}{0.58,0,0.82}
\usepackage{enumitem}
\setlist[itemize]{leftmargin=*}
\acmYear{2021}

\begin{document}

\title{\method: Multiresolution Tensor Completion from Partial and Coarse Observations}

\author{Chaoqi Yang$^1$, Navjot Singh$^1$, Cao Xiao$^{2}$, Cheng Qian$^{3}$, Edgar Solomonik$^{1}$, Jimeng Sun$^{1}$}
\affiliation{
  \institution{$^1$Department of Computer Science, University of Illinois Urbana-Champaign,  Urbana, IL 61801\\
$^2$Amplitude, San Francisco, CA 94105\quad
$^3$IQVIA, Cambridge, MA 02142}
\institution{$^1$\{chaoqiy2, navjot2, solomon2, jimeng\}@illinois.edu, $^2$danica.xiao@amplitude.com, $^3$alextoqc@gmail.com}
}


\begin{abstract}
    Existing tensor completion formulation mostly relies on partial observations from a single tensor. 
    However, tensors extracted from real-world data are often more complex due to: (i) {\bf Partial observation:} Only a small subset (e.g., 5\%) of tensor elements are available. 
    (ii) {\bf Coarse observation:} Some tensor modes only present coarse and aggregated patterns (e.g., monthly summary instead of daily reports). In this paper,
    we are given a subset of the tensor and some aggregated/coarse observations (along one or more modes) and seek to recover the original fine-granular tensor with low-rank factorization. We formulate a coupled tensor completion problem and propose an efficient \underline{M}ulti-resolution \underline{T}ensor \underline{C}ompletion (\method) model to solve the problem. Our \method model explores tensor mode properties and leverages the hierarchy of resolutions to recursively initialize an optimization setup, and optimizes on the coupled system using alternating least squares.
    \method ensures low computational and space complexity. We evaluate our model on two COVID-19 related spatio-temporal tensors. The experiments show that \method could provide 65.20\% and 75.79\% percentage of fitness (PoF) in tensor completion with only 5\% fine granular observations, which is 27.96\% relative improvement over the best baseline. To evaluate the learned low-rank factors, we also design a {tensor prediction task} for daily and cumulative disease case predictions, where \method achieves 50\% in PoF and 30\% relative improvements over the best baseline.
\end{abstract}

\maketitle

\section{Introduction}

Tensor completion is about estimating missing elements of multi-dimensional higher-order data (the dimensions of tensors are usually called ``mode''s). The CANDECOMP/PARAFAC (CP) and Tucker factorization approaches are commonly used for this data denois- ing and completion task, with applications ranging from image restoration~\cite{liu2012tensor,xu2013block} to healthcare data completion~\cite{acar2011,wang2015rubik},  recommendation systems~\cite{song2017multi}, link prediction~\cite{lacroix2018}, image fusion \cite{borsoi2020coupled} and spatio-temporal prediction~\cite{de2017,kargas2020stelar}.

In many real-world tensor data, missing values and observed data have unique characteristics:

\noindent
{\bf Partial observations:} Real-world tensor data has many missing values (only a small subset of the elements are available). This is a common assumption in most tensor completion works \cite{liu2012tensor,acar2011}. Many existing works estimate the missing elements purely based on data statistics \cite{hong2020generalized} by low-rank CP or Tucker model, while some also consider auxiliary domain knowledge \cite{lacroix2018,wang2015rubik}. 

\noindent
{\bf Coarse observations:} In many applications, tensor data is available at multiple granularity. For example, a location mode can be at zipcode, county, or state level, while a time mode can be at daily, weekly, or monthly level. Coarse-granular tensors are often fully available, while fine-granular tensors are usually incomplete with many missing elements. Few works \cite{almutairi2019prema,almutairi2020tendi} leverage the existing coarse level information to enhance tensor completion. 

\smallskip
\noindent{\bf Motivating example:} 
Given (i) {\bf partial observations:} a small subset (e.g., 5\%) of a fine-granular {\em location by disease by time} tensor, where each tensor element is the disease counts (by ICD-10 code\footnote{A disease can be represented with different coding standards: For example, International Classification of Disease version-10 (ICD-10) is a fine-granular code with over 100K dimensions. And Clinical Classification Software  for diagnosis codes (CCS) is a coarse-granular code with O(100) dimensions.}) in a geographical location (by anonymous location identifier) at a particular time (by date); (ii) {\bf coarse observations:} two coarse-granular tensors at {\em state by ICD-10 by date} level and {\em location by CCS by date} level, respectively. The problem is how to recover the fine-granular {\em location identifier by ICD-10 by date} tensor based on this partial and coarse information, where the mode aggregation mechanism is sometimes unknown, especially in healthcare \cite{park2014ludia}.


To capture the data characteristics, we identify the following technical challenges:
\begin{itemize}[leftmargin=*]
    \item {\bf Challenge in partial and coarse data fusion.} While heterogeneous information could compensate with each other and serve as potential supplements, it is challenging to combine different sources in a compact tensor factorization objective.
    \item {\bf Challenge in computational efficiency and accuracy}. With the increasing data volume and multidimensional structures, it becomes challenging to choose an accurate model initialization and a fast optimization kernel.
\end{itemize}
 

To address these challenges, we propose a Multiresolution Tensor Completion (\method) method with the following contributions.
\begin{itemize}[leftmargin=*]
    \item \method fuses both partial and coarse data into Frobenius norm based formulation, which is computationally efficient to deal with as they present a generic coupled-ALS form. 
    \item \method enables an effective initialization for model optimization to improve model efficiency and accuracy. We pack multiple linear equations into one joint normal equation during optimization and propose an ALS-based solver.
\end{itemize}
We evaluate \method on a spatio-temporal disease dataset and a public COVID keyword searching dataset. Extensive experiments demonstrate that for tensors with only 5\% observed elements, \method could provide 65.20\% and 75.79\% percentage of fitness (PoF) in the completion task, which shows relative improvement over the best baseline by 27.96\%. We also design a tensor prediction task to evaluate the learned factors, where our model reaches nearly 50\% PoF in disease prediction and outperforms the baselines by 30\%.

\section{Related Work on Tensors}
\subsection{Low-rank Tensor Completion}
Tensor completion (or tensor imputation) benefits various applications, which is usually considered a byproduct when dealing with missing data during the decomposition \cite{song2019tensor}. The low-rank hypothesis \cite{larsen2020practical} is often introduced to address the under-determined problem. Existing tensor completion works are mostly based on its partial \cite{acar2011,liu2012tensor} or coarsely \cite{almutairi2020tendi,almutairi2019prema} aggregated 
data. Many works have employed augmented loss functions based on statistical models, such as Wasserstein distance \cite{afshar2020swift}, missing values \cite{huang2020unified}, noise statistics (e.g., distribution, mean or max value) \cite{barak2016noisy} or the aggregation patterns \cite{yin2018joint}. Others have utilized other auxiliary domain knowledges, such as within-mode regularization \cite{kargas2020stelar,wang2015rubik}, $p$-norms \cite{lacroix2018}, pairwise constraints \cite{wang2015rubik}. Among which, a line of coupled matrix/tensor factorization \cite{almutairi2020tendi,de2017,bahargam2018constrained,huang2020unified,borsoi2020coupled} is related to our setting. One of the most noticeable works is coupled matrix and tensor factorization (CMTF) proposed by Acar. et al. \cite{acar2011all}, where a common factor is shared between a tensor and an auxiliary matrix. Another close work \cite{almutairi2019prema} recovered the tensor from only aggregated information. 

However, few works combine both the partial and coarse data by a coupled formulation. This paper fills in the gap and leverages both the partial and coarse information.

\subsection{Multi-scale and Randomized Methods}
With the rapid growth in data volume, multi-scale and randomized tensor methods become increasingly important for higher-order structures to boost efficiency and scalability  \cite{larsen2020practical,ozdemir2016multiscale}.

Multi-scale methods can interpolate from low-resolution to high-resolution \cite{schifanella2014multiresolution,sterck2013adaptive,park2020multiresolution} or operate first on one part of the tensor and then progressively generalize to the whole tensor \cite{song2017multi}. For example, Schifanella
 et al. \cite{schifanella2014multiresolution} exploits extra domain knowledge and develops a multiresolution method to improve CP and Tucker decomposition. Song et al. \cite{song2017multi} imputes the tensor from a small corner and gradually increases on all modes by multi-aspect streaming.
 
 The randomized methods are largely based on sampling, which accelerate the computation of over-determined least square problem \cite{larsen2020practical} in ALS for dense \cite{ailon2006approximate} and sparse \cite{eshragh2019lsar} tensors by effective strategies, such as Fast Johnson-Lindenstrauss Transform (FJLT) \cite{ailon2006approximate}, leverage-based sampling \cite{eshragh2019lsar}  and tensor sketching.
 
 Unlike previous methods, our model does not require a domain-specific mode hierarchy. Our multiresolution factorization relies on a heuristic mode sampling technique (mode is continuous or categorical), which is proposed to initialize the optimization phase such that fewer iterations would be needed in high resolution. Also, we consider a more challenging coupled factorization problem.

\section{Preliminary}
\subsection{Notations}
We use plain letters for scalar, such as $x$ or $X$, boldface uppercase letters for matrices, e.g., $\mathbf{X}$, and Euler script letter for tensors and sets, e.g., $\mathcal{X}$. Tensors are multidimensional arrays indexed by three or more indices (modes). For example, an $N$-th order tensor $\mathcal{X}$ is an $N$-dimensional array of size $I_1\times\cdots\times I_N$, where $\mathcal{X}({i_1,\dots,i_N})$ is the element at the $(i_1,\cdots,i_N)$-th position. For matrices, $\mathbf{X}({:,r})$ and $\mathbf{X}({r,:})$ are the $r$-th column and row, and $\mathbf{X}({i,j})$ is for the $(i,j)$-th element. Indices typically start from $1$, e.g., $\mathbf{X}(1,:)$ is the first row of the matrix. In this paper, we work with a three-order tensor for simplicity,
while our method could also be applied to higher-order tensors. 

\subsubsection{CANDECOMP/PARAFAC (CP) Decomposition.} One of the common compression methods for tensors is CP decomposition \cite{hitchcock1927expression,carroll1970analysis}, which approximates a tensor by multiple rank-one components. For example, let $\mathcal{X}\in\mathbb{R}^{I_1\times I_2\times I_3}$ be an arbitrary $3$-mode tensor of CP rank $R$, it can be expressed exactly by three factor matrices $\mathbf{U}\in\mathbb{R}^{I_1\times R},\mathbf{V}\in\mathbb{R}^{I_2\times R},\mathbf{W}\in\mathbb{R}^{I_3\times R}$ as
\begin{equation}
    \mathcal{X} = \sum_{r=1}\mathbf{U}(:,r)\circ\mathbf{V}(:,r)\circ \mathbf{W}(:,r) = \llbracket\mathbf{U},\mathbf{V},\mathbf{W}\rrbracket,
\end{equation}
where $\circ$ is the vector outer product, and the double bracket notation $\llbracket\,\cdot,\cdot,\cdot\rrbracket$ is the Kruskal operator.
 
\subsubsection{Mode Product.} For the same tensor $\mathcal{X}$, a coarse observation is given by conducting a tensor mode product with an aggregation matrix (for example, on the first mode), $\mathbf{P}_1\in\mathbb{R}^{J_1\times I_1}$, where $J_1<I_1$,
\begin{align}
    \mathcal{C}^{(1)} &= \mathcal{X}\times_1\mathbf{P}_1 =  \sum_{r=1}\mathbf{P}_1\mathbf{U}(:,r)\circ\mathbf{V}(:,r)\circ \mathbf{W}(:,r) \notag\\
    &= \llbracket\mathbf{P}_1\mathbf{U},\mathbf{V},\mathbf{W}\rrbracket\in\mathbb{R}^{J_1\times I_2\times I_3},
\end{align}
where we use superscripts, e.g., $\mathcal{C}^{(1)}$, to denote the aggregated tensor. Other tensor background may be found in Appendix. We organize the symbols in Table~\ref{tb:notations}.

\begin{table}[h!] \small \caption{Notations used in \method}\vspace{-3mm}
\centering
\resizebox{0.49\textwidth}{!}{
\begin{tabular}{l|l} \toprule \textbf{Symbols} & \textbf{Descriptions} \\ 
	\midrule 
	$R$ & tensor CP rank \\
	$\odot$, $*$ & tensor Khatri-Rao/Hadamard product \\
	$\llbracket\cdot,\cdot,\cdot\rrbracket$ & tensor Kruskal operator \\
	$\mathcal{J}$ & the set of coarse tensors \\
	$\mathcal{P}^*,\tilde{\mathcal{P}},\mathcal{P}$ & known/unknown/full aggregation matrix set \\
	$\mathbf{P}_1\in\mathbb{R}^{J_1\times I_1}, \mathbf{P}_2\in\mathbb{R}^{J_2\times I_2}$ & aggregation matrices \\
	$\mathcal{X},\tilde{\mathcal{X}},\mathcal{R}\in\mathbb{R}^{I_1\times I_2\times I_3}$ & the original/interim/reconstructed tensor \\
	$\mathcal{C}^{(1)}\in\mathbb{R}^{J_1\times I_2\times I_3}$, $\mathcal{C}^{(2)}\in\mathbb{R}^{I_1\times J_2\times I_3}$ & two coarse observation tensors \\
	$\tilde{\mathbf{X}}_i,\mathbf{C}_i^{(1)},\mathbf{C}_i^{(2)},\mathbf{M}_i,\mathbf{R}_i,~i=1,2,3$ & tensor unfoldings/matricizations\\
	$\mathcal{M}\in\{0,1\}^{I_1\times I_2\times I_3}$ & tensor mask \\
	$\mathbf{U}\in\mathbb{R}^{I_1\times R}, \mathbf{V}\in\mathbb{R}^{I_2\times R}, \mathbf{W}\in\mathbb{R}^{I_3\times R}$ & factor matrices \\
	$\mathbf{Q}_1\in\mathbb{R}^{J_1\times R}, \mathbf{Q}_2\in\mathbb{R}^{J_2\times R}$ & auxiliary factor matrices \\
	$(\cdot)_{(d_t)},~\mathcal{S}_{(d_{t}, ~\cdot)}$ & entity or index set at $d_t$ resolution \\
	$\tilde{\mathbf{U}},\tilde{\mathbf{V}},\tilde{\mathbf{W}},\tilde{\mathbf{Q}}_1,\tilde{\mathbf{Q}}_2$ & initializations for the factors \\
    \bottomrule 
    \end{tabular}} \vspace{-3mm}
	\label{tb:notations} \end{table}

\subsection{Problem Definition}

\begin{problem}[\textbf{Tensor Completion with Partial and Coarse Observations}]
For an unknown tensor $\mathcal{X}\in\mathbb{R}^{I_1\times I_2 \times I_3}$, we are given
\begin{itemize}
    \item a partial observation $\mathcal{M}*\mathcal{X}$, parameterized by mask $\mathcal{M}\in\{0,1\}^{I_1\times I_2 \times I_3}$;
    \item a set of known aggregation matrices $\mathcal{P}^* \subseteq \mathcal{P}=\{\mathbf{P}_s: {\mathbf{P}_s\in\{0,1\}^{J_s \times I_s}},\\ ~J_s < I_s, ~s=1,2,3\}$, each column of $\mathbf{P}_s$ is one-hot;
    \item a set of coarse tensors $\mathcal{J}=\{\mathcal{C}^{\mathcal{H}}:\mathcal{C}^{\mathcal{H}}=\mathcal{X}\times_{h\in\mathcal{H}}\mathbf{P}_h,~\mathcal{H}\in\mbox{Powerset}\{1,2,3\}\setminus \varnothing, ~\mathbf{P}_h\in\mathcal{P}\}$. For example, an coarse tensor aggregated on both the first and the second mode can be written as $\mathcal{C}^{(1,2)}=\mathcal{X}\times_{h\in\{1,2\}}\mathbf{P}_h=\mathcal{X}\times_1\mathbf{P}_1\times_2\mathbf{P}_2$.

\end{itemize} 
The problem is to find low-rank CP factors $\mathbf{U}, \mathbf{V}, \mathbf{W}$, 
such that the following (Frobenius norm based) loss is minimized,
\begin{align}
    \mathcal{L}(\mathcal{M}*\mathcal{X},\mathcal{J},\mathcal{P}^*;~\mathbf{U,V,W}, {\tilde{\mathcal{P}}})=~ \left\|\mathcal{M} * (\mathcal{X} - \llbracket \mathbf{U,V,W} \rrbracket)\right\|_F^2& \notag\\
    +~\sum_{\mathcal{H}}\lambda_{\mathcal{H}}\|\mathcal{C}^\mathcal{H}-\llbracket\mathbf{U},\mathbf{V},\mathbf{W}\rrbracket\times_{h\in\mathcal{H}}\mathbf{P}_h\|_F^2&,
\end{align}
where $\tilde{\mathcal{P}} = \mathcal{P}\setminus \mathcal{P}^*$ is the set of unknown aggregation matrices, $\mathcal{H}$ enumerates the index of $\mathcal{J}$, and $\lambda_\mathcal{H}$ contains the weights that specify the importance of each coarse tensor. We separate the knowns and unknowns (to estimate) by semicolon ``;''. Note that, the same mode cannot be aggregated in all coarse tensors.

The problem reduces to the conventional CP tensor completion problem if the set $\mathcal{J}$ of coarse tensors is empty, which has been explored extensively in the previous literature \cite{liu2012tensor,acar2011,10.1016/j.parco.2015.10.002}.

\end{problem}

\subsection{A Motivating Example}
As an example, we use the aforementioned COVID disease tensor to motivate our application. 
For an unknown {\em location identifier by ICD-10 by date} tensor $\mathcal{X}\in\mathbb{R}^{I_1\times I_2 \times I_3}$, we are given
\begin{itemize}
    \item a partial observation $\mathcal{M}*\mathcal{X}$, parametrized by mask $\mathcal{M}$; 
    \item a coarse {\em state by ICD-10 by date} tensor, $\mathcal{C}^{(1)}\in\mathbb{R}^{J_1\times I_2\times I_3}$, where $J_1<I_1$ (satisfying $\mathcal{C}^{(1)}=\mathcal{X}\times_1\mathbf{P}_1$ with  an unknown $\mathbf{P}_1\in\mathbb{R}^{J_1\times I_1}$, since we have no knowledge on the mapping from anonymized location identifier to state);
    \item a coarse {\em location identifier by CCS by date} tensor, $\mathcal{C}^{(2)}\in\mathbb{R}^{I_1\times J_2\times I_3}$, where $J_2<I_2$ and it satisfies $\mathcal{C}^{(2)}=\mathcal{X}\times_2\mathbf{P}_2$ with known {\em ICD-10 to CCS} mapping $\mathbf{P}_2\in\mathbb{R}^{J_2\times I_2}$.
\end{itemize} 
We seek to find low-rank CP factors $\mathbf{U}, \mathbf{V}, \mathbf{W}$, such that the loss is minimized over parameters ($\mathbf{U,V,W,P_1}$),
\begin{align} \label{eq:original_problem}
    \mathcal{L}(\mathcal{M}*\mathcal{X},\mathcal{C}^{(1)},\mathcal{C}^{(2)},\mathbf{P}_2;\mathbf{U,V},&\mathbf{W},\mathbf{P}_1)= \left\|\mathcal{M} * (\mathcal{X} - \llbracket \mathbf{U,V,W} \rrbracket)\right\|_F^2 \notag\\
    &+ ~\lambda_1\left\|\mathcal{C}^{(1)} - \llbracket \mathbf{P}_1 \mathbf{U,V,W} \rrbracket\right\|_F^2& \notag\\
    &+~\lambda_2\left\|\mathcal{C}^{(2)} - \llbracket \mathbf{U},\mathbf{P}_2\mathbf{V,W} \rrbracket\right\|_F^2,
\end{align}
\noindent w.l.o.g., we use this example to introduce our methodology.


\subsection{Solution Outline}


\subsubsection{Coupled Tensor Decomposition.}


To address the above objective, we consider to successively approximate the low-rank completion problem by coupled tensor decomposition. We first introduce two transformations and reformulate the objective into a coupled decomposition form:
\begin{itemize}
    \item {\bf Parameter replacement.} We replace $\mathbf{Q}_1=\mathbf{P}_1\mathbf{U}$ and $\mathbf{Q}_2=\mathbf{P}_2\mathbf{V}$ as two new auxiliary factor matrices\, and add extra constrains later since $\mathbf{Q}_2$ and $\mathbf{V}$ are related by a known $\mathbf{P}_2$ (but $\mathbf{P}_1$ is unknown, we does not consider the constraint for $\mathbf{Q}_1$ and $\mathbf{U}$).
    \item {\bf Interim Tensor.} For the first term in Equation~\eqref{eq:original_problem}, we construct an interim tensor $\tilde{\mathcal{X}}$ from an expectation-maximization (EM) \cite{acar2011} approach,
{ \begin{align} 
    &\mathcal{M} * (\mathcal{X} - \llbracket \mathbf{U,V,W}\rrbracket) \notag\\ 
    = & ~\mathcal{M} * \mathcal{X} + (1-\mathcal{M}) * \llbracket \mathbf{U,V,W} \rrbracket - \llbracket \mathbf{U,V,W} \rrbracket \notag\\
    \approx & \underbrace{\left[\mathcal{M} * \mathcal{X} + (1-\mathcal{M}) * \llbracket \mathbf{U}^{k},\mathbf{V}^{k},\mathbf{W}^{k} \rrbracket\right]}_{\mbox{the interim tensor $\tilde{\mathcal{X}}$ for next iteration}} -~ \llbracket  \mathbf{U},\mathbf{V},\mathbf{W} \rrbracket,
\end{align}}

where $\mathbf{U}^{k},\mathbf{V}^{k},\mathbf{W}^{k}$ are the learned factors at the $k$-th iteration. The reconstruction $\mathcal{R}=\llbracket \mathbf{U}^{k},\mathbf{V}^{k},\mathbf{W}^{k} \rrbracket$ is expensive for large tensors. However, we do not need to compute it explicitly. In Section~\ref{sec:optimization}, we introduce a way to use the implicit form.
\end{itemize}

 \noindent Afterwards, the objective in Equation~\eqref{eq:original_problem} is approximated by a coupled tensor decomposition form (at the $(k+1)$-th iteration),
\begin{align}
    \mathcal{L}(\tilde{\mathcal{X}},\mathcal{C}^{(1)},\mathcal{C}^{(2)},\mathbf{P_2};\mathbf{U,V,W},\mathbf{Q}_1,\mathbf{Q}_2)=~\left\| \tilde{\mathcal{X}} - \llbracket \mathbf{U,V,W} \rrbracket\right\|_F^2&\notag\\
    +\lambda_1\left\|\mathcal{C}^{(1)} - \llbracket \mathbf{Q}_1, \mathbf{V,W} \rrbracket\right\|_F^2& \notag \\
    +\lambda_2\left\|\mathcal{C}^{(2)} - \llbracket \mathbf{U},\mathbf{Q}_2,\mathbf{W} \rrbracket\right\|_F^2& \label{eq:original},
\end{align}
with constraint $\mathbf{Q_2}=\mathbf{P_2V}$ and 
the interim tensor $\tilde{\mathcal{X}}$,
\begin{equation} \label{eq:implicit_form}
       \tilde{\mathcal{X}} =\mathcal{M} * \mathcal{X} + (1-\mathcal{M}) * \llbracket \mathbf{U}^{k},\mathbf{V}^{k},\mathbf{W}^{k} \rrbracket. \\
\end{equation}
Now, our goal turns into optimizing a coupled decomposition problem with five shared factors, among which $\mathbf{U,V,W}$ are our focus.

\subsubsection{Optimization Idea.}
\begin{figure}
    \centering
    \includegraphics[width=0.49\textwidth]{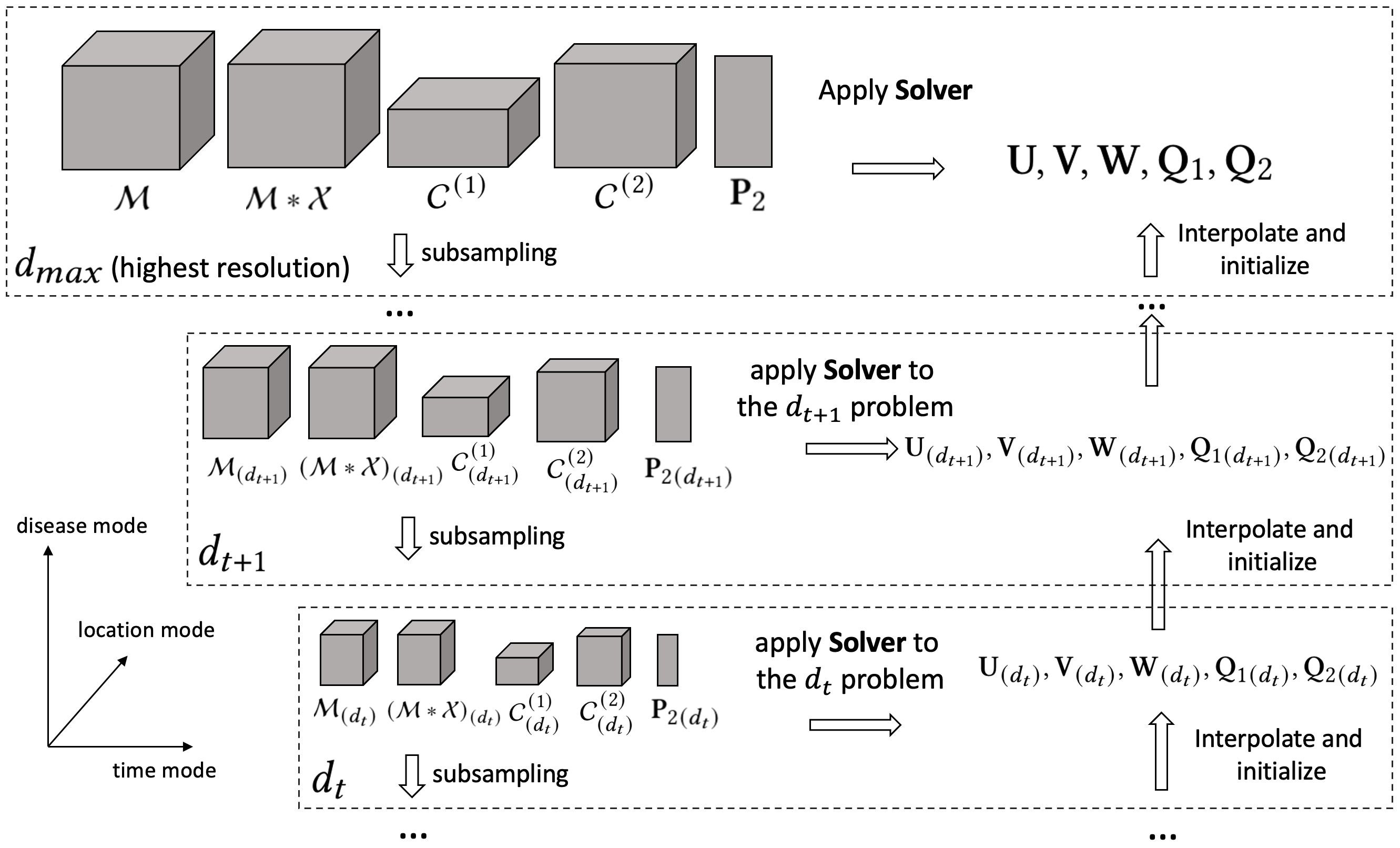}
    \caption{ \method Framework {\normalfont follows a multiresolution recursive flow. First, we apply subsampling on all accessible information (i.e., the tensors and the known aggregation matrices) into the lowest resolution. Then, we solve the low resolution problem by applying the optimization solver. Next, we interpolate the solution into the higher resolution to initialize the high resolution factors. We repeat this process and find  a good initialization for the original fine-granular problem. The mode {\bf subsampling and interpolation} are discussed in Section~\ref{sec:multiresolution}; the {\bf solver} is developed in Section~\ref{sec:optimization}}.}
    \label{fig:framework}
\end{figure}

As a solution, we adopt a multiresolution strategy, as in Figure~\ref{fig:framework}, which leverages a hierarchy of resolution to perform the coupled optimization effectively.
The {\bf multiresolution factorization} part (in Section~\ref{sec:multiresolution}) serves as the outer loop of our solution. At each resolution, we derive a joint normal equation for each factor and develop an efficient ALS-based {\bf solver} (in Section~\ref{sec:optimization}) to handle multiple linear systems simultaneously.




\section{Multiresolution Factorization} \label{sec:multiresolution}
Let us denote the resolution hierarchy as $\{d_{max},\dots,d_{t+1},d_t,\dots\}$, where $d_{max}$ is the highest resolution (i.e., original tensors). The multiple resolutions are computed via recursive subsampling along three tensor modes simultaneously. The underlying ideas of the multiresolution factorization are: (i) to subsample high-resolution information into low resolution and solve a low-resolution problem; (ii) to then interpolate the low-resolution solution factors for initializing the high-resolution problem, such that fewer iterations are needed in high resolution.

\subsection{High ($d_{t+1}$) to Low ($d_t$) Resolution} 

 Assume that $d_{t+1}$ is a generic notation for high resolution and $d_{t}$ is for the neighboring low resolution, we use {\em subscript-parenthesis} notation, e.g., $_{(d_t)}$ to denote the entities at each resolution.

{\bf Objective at $d_{t+1}$.} The goal is to find optimal factor matrices:  $\mathbf{U}_{(d_{t+1})}\in\mathbb{R}^{I_{1(d_{t+1})}\times R}$, $\mathbf{V}_{(d_{t+1})}\in\mathbb{R}^{I_{2(d_{t+1})}\times R}$, $\mathbf{W}_{(d_{t+1})}\in\mathbb{R}^{I_{3(d_{t+1})}\times R}$, $\mathbf{Q}_{1(d_{t+1})}\in\mathbb{R}^{J_{1(d_{t+1})}\times R}$, $\mathbf{Q}_{2(d_{t+1})}\in\mathbb{R}^{J_{2(d_{t+1})}\times R}$, such that they minimize the loss at $d_{t+1}$ (the form could refer to Equation~\eqref{eq:original}),
{\begin{align}
    \mathcal{L}_{(d_{t+1})} = \mathcal{L}\left(\tilde{\mathcal{X}}_{(d_{t+1})},\mathcal{C}^{(1)}_{(d_{t+1})},\mathcal{C}^{(2)}_{(d_{t+1})},\mathbf{P}_{2(d_{t+1})};\mathbf{U}_{(d_{t+1})}, \mathbf{V}_{(d_{t+1})},\right.\notag\\
    \left.\mathbf{W}_{(d_{t+1})},\mathbf{Q}_{1(d_{t+1})},\mathbf{Q}_{2(d_{t+1})}\right).\notag
\end{align}}
In the following, we conceptulize the solution to the objective.
\begin{itemize}
    \item[(1)] {\bf Subsampling}: We first subsample all accessible information (tensors and the known aggregation matrices) from $d_{t+1}$ to $d_{t}$ at once,
    \begin{align}
    \mathcal{M}_{(d_{t+1})}, {(\mathcal{M}*\mathcal{X})}_{(d_{t+1})},\mathcal{C}^{(1)}_{(d_{t+1})},\mathcal{C}^{(2)}_{(d_{t+1})},\mathbf{P}_{2(d_{t+1})} ~\stackrel{subsample}{\longrightarrow}&\notag \\
    \mathcal{M}_{(d_t)},{(\mathcal{M}*\mathcal{X})}_{(d_{t})},\mathcal{C}^{(1)}_{(d_{t})},\mathcal{C}^{(2)}_{(d_{t})},\mathbf{P}_{2(d_{t})}&.
    \end{align}
    \item[(2)] {\bf Solve the objective at $d_{t}$}: Next, we solve the $d_t$-problem recursively and obtain smaller factor matrices ($\mathbf{U}_{(d_{t})}\in\mathbb{R}^{I_{1(d_{t})}\times R}$, $\mathbf{V}_{(d_{t})}\in\mathbb{R}^{I_{2(d_{t})}\times R}$, $\mathbf{W}_{(d_{t})}\in\mathbb{R}^{I_{3(d_{t})}\times R}$,
$\mathbf{Q}_{1(d_{t})}\in\mathbb{R}^{J_{1(d_{t})}\times R}$,$\mathbf{Q}_{2(d_{t})}\in\mathbb{R}^{J_{2(d_{t})}\times R}$), where the $d_t$-objective can be similarly written as
\begin{align} \label{eq:dt_problem}
    \mathcal{L}_{(d_t)}
    &=~ \left\| \tilde{\mathcal{X}}_{(d_t)} - \llbracket \mathbf{U}_{(d_t)},\mathbf{V}_{(d_t)},\mathbf{W}_{(d_t)} \rrbracket\right\|_F^2 \notag\\
    &+\lambda_1\left\|\mathcal{C}^{(1)}_{(d_t)} - \llbracket \mathbf{Q}_{1(d_t)}, \mathbf{V}_{(d_t)},\mathbf{W}_{(d_t)} \rrbracket\right\|_F^2 \notag \\
    &+\lambda_2\left\|\mathcal{C}^{(2)}_{(d_t)} - \llbracket \mathbf{U}_{(d_t)},\mathbf{Q}_{2(d_t)},\mathbf{W}_{(d_t)} \rrbracket\right\|_F^2,
\end{align}
the interim tensor $\tilde{\mathcal{X}}_{(d_t)}$ is calculated by
\begin{equation}
        \tilde{\mathcal{X}}_{(d_t)} =(\mathcal{M}* \mathcal{X})_{(d_t)} + (1-\mathcal{M}_{(d_t)}) * \llbracket \mathbf{U}_{(d_t)}^{k},\mathbf{V}_{(d_t)}^{k},\mathbf{W}_{(d_t)}^{k} \rrbracket,
\end{equation}
and the $d_t$ resolution constraint is $\mathbf{Q}_{2(d_t)}=\mathbf{P}_{2(d_t)}\mathbf{V}_{(d_t)}$.

\smallskip
\item[(3)] {\bf Interpolation:} Then, we interpolate the $d_t$ solutions to initialize the $d_{t+1}$ problem. 
More precisely, four factor matrices are obtained by interpolation:
\begin{align}
    \mathbf{U}_{(d_{t})},\mathbf{V}_{(d_{t})},\mathbf{W}_{(d_{t})},\mathbf{Q}_{1(d_{t})} ~\stackrel{interpolate}{\longrightarrow}& \\ \tilde{\mathbf{U}}_{(d_{t+1})},\tilde{\mathbf{V}}_{(d_{t+1})},\tilde{\mathbf{W}}_{(d_{t+1})},\tilde{\mathbf{Q}}_{1(d_{t+1})}&.
\end{align}

And one is calculated by the known aggregation matrix according to the $d_{t+1}$ constraint:
\begin{equation}
    \tilde{\mathbf{Q}}_{2(d_{t+1})} = \mathbf{P}_{2(d_{t+1})} \tilde{\mathbf{V}}_{(d_{t+1})}.
\end{equation}

\item[(4)] {\bf Optimization at $d_{t+1}$}: With the initialized factor matrices, we finally can apply our optimization algorithm recursively on the $d_{t+1}$ objective to improve the factors.
\end{itemize}

These four steps are the key to the paper. The multiresolution factorization is conducted by recursively applying four steps within Step (2). In this section, Step (1) and (3) are of our focus, while we will specify Step (4) in Section 5. The subsampling and interpolation are closely related, and we elucidate them below. 

\subsection{Resolution Transition from $d_{t+1}$ to $d_t$}
 In our example, three modes are presented with $d_{t+1}$-specific mode sizes: (i) location mode: anonymized location identifier $I_{1(d_{t+1})}$, state $J_{1(d_{t+1})}$; (ii) disease mode: ICD-10 $I_{2(d_{t+1})}$, CCS $J_{2(d_{t+1})}$; (iii) time mode: date $I_{3(d_{t+1})}$. To make it clear, we use the term {\bf aspect} to indicate a tensor mode or its projection/aggregation at a coarse granularity. 
 
 In total, we have three modes in our problem and five different aspects, corresponding to five factors. We call a mode with two aspects an {\em aggregated mode} and a mode with one aspect as an {\em single mode}. For example, the location and disease modes are aggregated modes, and the time mode is a single mode. Because location mode can be at zipcode or state aspect; disease mode can be at ICD-10 aspect (fine-granular) or CCS aspect (coarse-granular). But time mode is only represented at the date aspect.  
 
{\em Mode subsampling is actually conducted on each aspect.} Tensors and aggregation matrices are subsampled altogether in the same way when they share the same aspect. For example, on subsampling ICD-10 aspect, $\mathcal{M}$ (the second mode), $\mathcal{M}*\mathcal{X}$ (the second mode), $\mathcal{C}^{(1)}$ (the second mode) and $\mathbf{P}_2$ (the second dimension) are targeted. We consider different subsampling scenarios.
\begin{itemize}
    \item {\bf Continuous or categorical.} For single modes, we provide two strategies based on mode property: (i) mode with continuous or smooth information (called continuous mode), e.g., date; (ii) mode with categorical information (called categorical mode).
    \item {\bf Known or unknown aggregation.} For aggregated mode, we provide different treatments depending on whether the associated aggregation matrix $\mathbf{P}$ is (i) known or (ii) unknown.
\end{itemize}

In fact, the strategies for aggregated modes are based on the single mode (continuous or categorical) strategy. We first introduce the strategies on a single mode.

\subsection{Single Continuous Mode}\label{sec:continuous} The transition for this mode type is under the assumption of smoothness: information along the mode contains some spatial or {temporal} consistency. For a continuous mode, we subsample the $d_{t+1}$ index set (of the aspect) to $d_t$ and later use neighborhood smoothing to interpolate the $d_t$  factors back to initialize $d_{t+1}$ factors. 

\begin{figure}
    \centering
    \includegraphics[width=0.47\textwidth]{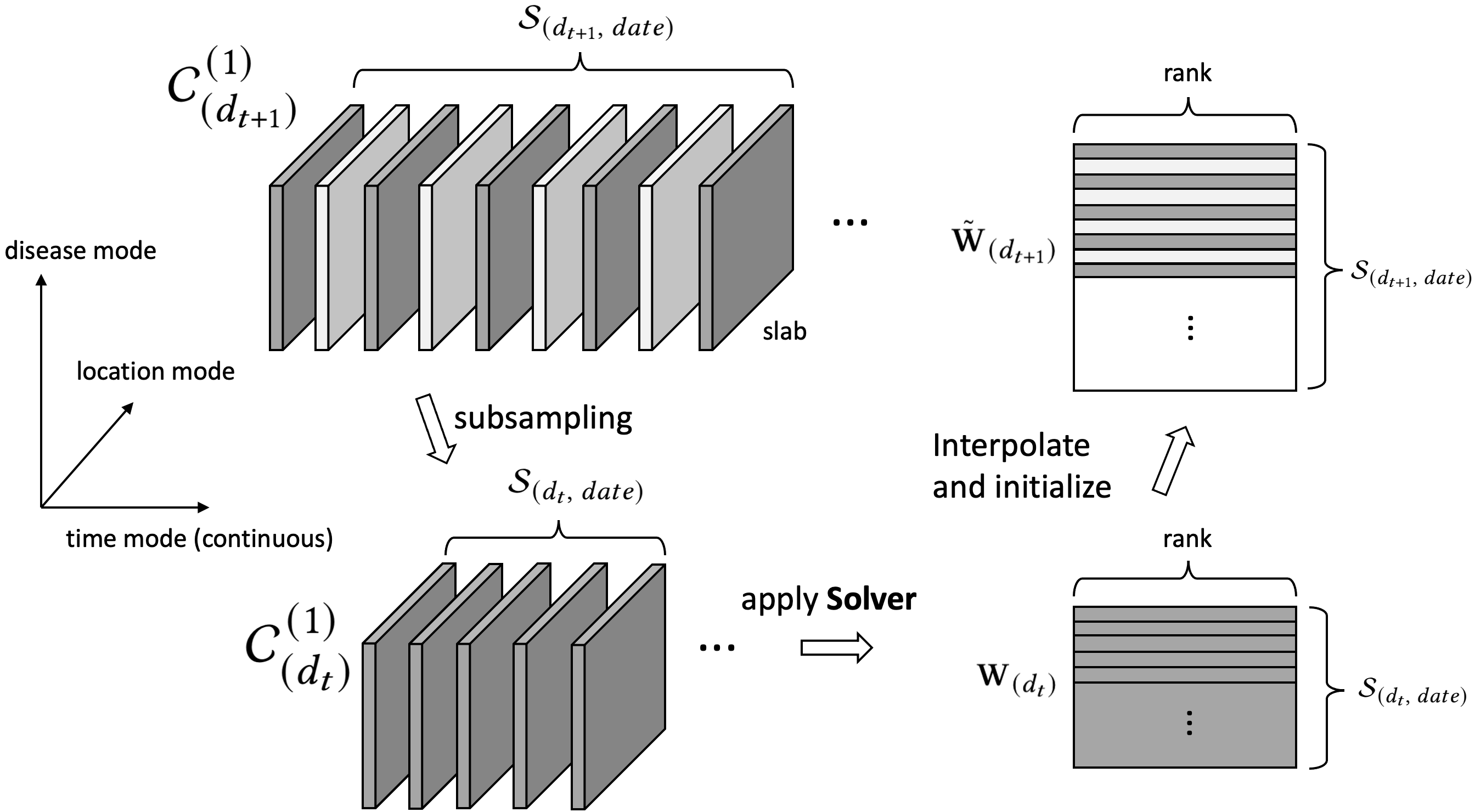}
    \caption{Continuous Mode. {\normalfont We use the time mode (corresponding to factor $\mathbf{W}$) and one of the involved tensors, $\mathcal{C}^{(1)}$, as an example. Here, we fix other modes and only discuss the transition along the time mode (in real practice, we do subsampling along three modes simultaneously).}}
    \label{fig:continuous_illustration}
\end{figure}

For example, the time mode assumes smoothness, since within a time period, the disease statistics over the same geographical regions will not change abruptly. See the diagram in Figure~\ref{fig:continuous_illustration}.

\subsubsection{Subsampling.} We define
\begin{equation}
    \mathcal{S}_{(d_{t+1}, ~date)} = \{1,2,\dots, I_{3(d_{t+1})}\}
\end{equation}
to be the index set of the date aspect at $d_{t+1}$-resolution. 
We sample at a regular interval as
\begin{equation}
    \mathcal{S}_{(d_t,~date)} = \{1,3,5,\dots\},
\end{equation}
yielding a $d_t$-resolution index set of size $\lceil\frac{I_{3(d_{t+1})}}{2}\rceil$. 

For example, $\mathcal{C}^{(1)}_{(d_{t+1})}$ includes the date aspect,  thus it will be subsampled using the index set $\mathcal{S}_{(d_t,~date)}$.

\subsubsection{Interpolation.} Factor $\mathbf{W}_{(d_t)}\in\mathbb{R}^{I_{3(d_t)}\times R}$ corresponds to the date aspect at $d_t$. The initialization of $d_{t+1}$ factor $\tilde{\mathbf{W}}_{(d_{t+1})}\in\mathbb{R}^{I_{3(d_{t+1})}\times R}$ is given by
\begin{align}
    \tilde{\mathbf{W}}_{(d_{t+1})}(2j,:) &= \frac12\left(\mathbf{W}_{(d_t)}(j,:)+\mathbf{W}_{(d_t)}(j+1,:)\right) \label{eq:interpolate1},\\
    \tilde{\mathbf{W}}_{(d_{t+1})}(2j+1,:) &= \mathbf{W}_{(d_t)}(j,:),\quad 0\leq j\leq\left\lfloor\frac{I_{3(d_{t+1})}-1}{2}\right\rfloor, \label{eq:interpolate2}
\end{align}
if $I_{3(d_{t+1})}$ is an even number, then
\begin{equation} \label{eq:interpolate3}
    \tilde{\mathbf{W}}_{(d_{t+1})}(I_{3(d_{t+1})},:) = \mathbf{W}_{(d_t)}\left(\frac{I_{3(d_{t+1})}}{2},:\right).
\end{equation}

\subsubsection{Intuition.} As shown in Figure~\ref{fig:continuous_illustration}, from $\mathcal{C}^{(1)}_{(d_{t+1})}$ to $\mathcal{C}^{(1)}_{(d_{t})}$, only the ``dark'' slabs are selected for the $d_t$-resolution. Tensor slabs of $\mathcal{C}^{(1)}_{(d_{t})}$ along the date aspect are controlled by $\mathbf{W}_{(d_{t})}$, where one row corresponds to one slab. After solving the $d_t$ problem, the solution factor $\mathbf{W}_{(d_{t})}$ provides an estimate for those ``dark'' slabs. To  also estimate the ``light'' slabs, we assume that it could be approximated by the neighboring ``dark'' slabs. The actual interpolation follows Equation~\eqref{eq:interpolate1}, \eqref{eq:interpolate2}, \eqref{eq:interpolate3}, where we average neighboring rows.

\subsection{Single Categorical Mode} For a categorical mode, we also want to subsample the $d_{t+1}$ index set to $d_{t}$. However, we now focus on the density of the corresponding tensor slabs. Later on, we interpolate the $d_{t}$ factors back to $d_{t+1}$.  

\begin{figure}
    \centering
    \includegraphics[width=0.49\textwidth]{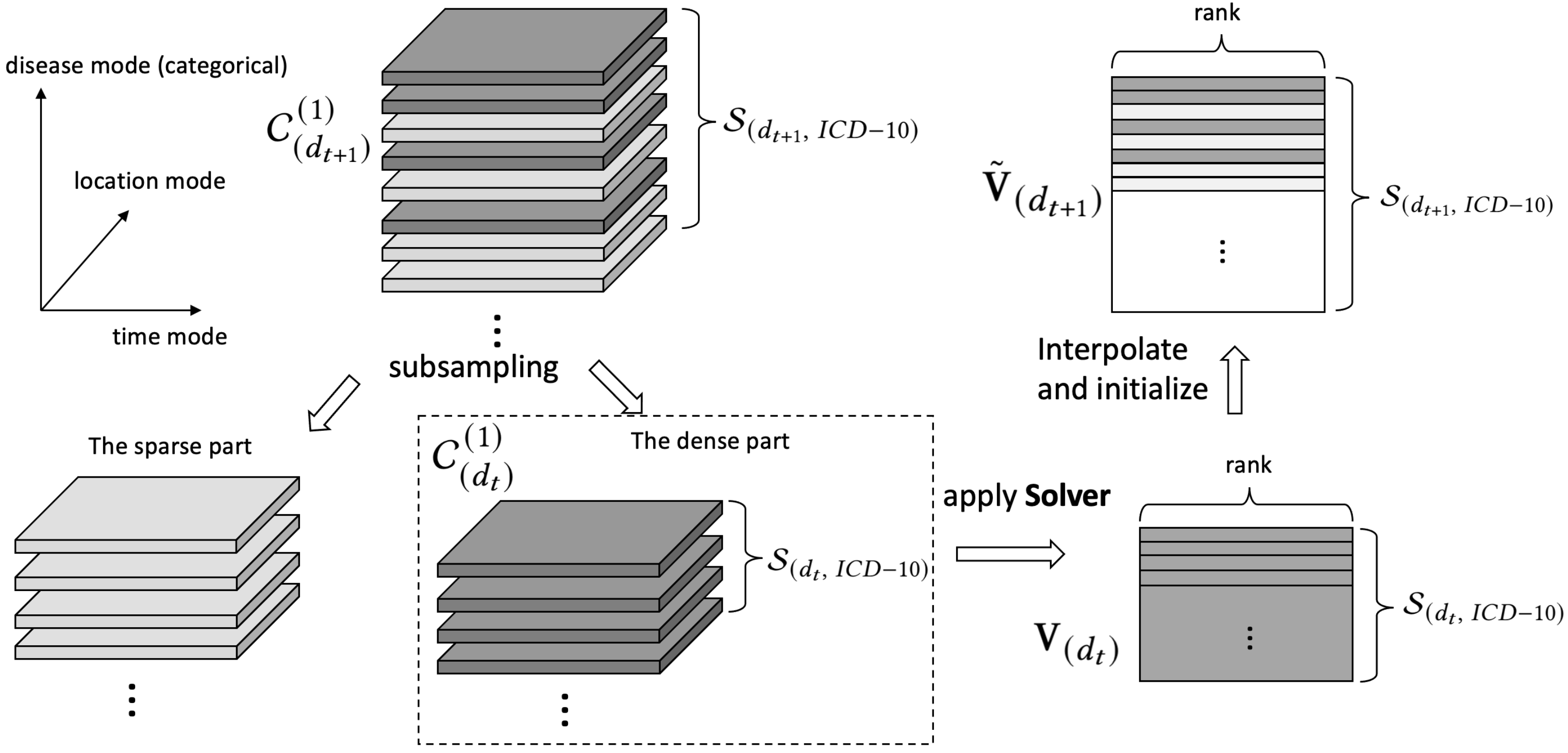}
    \caption{Categorical Mode. {\normalfont We use the ICD-10 aspect at disease mode (corresponding to factor $\mathbf{V}$) and one of the involved tensors, $\mathcal{C}^{(1)}$, as an example. Similarly, only one aspect is discussed here.}}
    \label{fig:categorical_illustration}
\end{figure}

We use the ICD-10 (fine-granular disease code) aspect as an example. 
This aspect corresponds to factor $\mathbf{V}_{(d_{t+1})}$ at $d_{t+1}$, for which we provide an illustration in Figure~\ref{fig:categorical_illustration}.


The ICD-10 information could be obtained from one coarse tensor, e.g., $\mathcal{C}^{(1)}_{(d+1)}$, at $d_{t+1}$. Along this aspect, we count the non-zero elements over each slab and obtain a count vector; each element stores the counts of one slab, corresponding to the ICD-10 code,
\begin{equation}
    Cnt_{(d_{t+1}, ~ICD-10)} = [Cnt_1,Cnt_2,\dots,Cnt_{I_{2(d_{t+1})}}].
\end{equation}

\subsubsection{Subsampling.} We define
\begin{equation}
    \mathcal{S}_{(d_{t+1},~ICD-10)} = \{1,2,3,\dots,I_{2(d_{t+1})}\}
\end{equation}
to be the index set of ICD-10 aspect at $d_{t+1}$-resolution. 

We order the count vector and find the indices corresponding to the largest half (composing the dense part of the tensor, while the unselected half is the sparse part), denoted as $\mathcal{S}_{(d_t,~ICD-10)}$, as the index set at $d_t$, so that
\begin{equation}
    Cnt_i \geq Cnt_j,\quad \forall i\in \mathcal{S}_{(d_t,~ICD-10)}, ~j\in \mathcal{S}_{(d_{t+1},~ICD-10)}\setminus \mathcal{S}_{(d_t,~ICD-10)}.
\end{equation}
The tensors and aggregation matrices that include ICD-10 aspect, e.g., $\mathcal{C}^{(1)}_{(d_{t+1})}$, are subsampled along this mode.

\subsubsection{Interpolation.} $\mathbf{V}_{(d_t)}\in\mathbb{R}^{{I_{2(d_t)}}\times R}$ is the corresponding solution factor for ICD-10 at $d_t$. We get the $d_{t+1}$ initialization $\tilde{\mathbf{V}}_{(d_{t+1})}\in\mathbb{R}^{{I_{2(d_{t+1})}}\times R}$ as follows. The rows of $\tilde{\mathbf{V}}_{(d_{t+1})}$ are direct copies from the rows of $\mathbf{V}_{(d_t)}$ if the {corresponding indices} are selected in $\mathcal{S}_{(d_t,~ICD-10)}$, otherwise, we fill the rows with random entries (in the paper, we use i.i.d. samples from $[-1,1]$, and users might customize their configuration based on the applications).

\subsubsection{Intuition.} As for the categorical mode (in Figure~\ref{fig:categorical_illustration}), the rationality of our strategy comes from two underlying heuristics: (i) the dense half of the indices is likely to account for a large volume of the tensor, such that after the optimization at $d_t$, we already recover a dominant part at $d_{t+1}$; (ii) the selected indices (slabs) preserve the dense half of the tensor, while from randomized algorithm perspective, {a slab may be regarded as a training sample}, composited by the other factor matrices. Preserving those dense slabs would be beneficial for the estimation of other factors.

\subsection{Aggregated Mode} 


We are now ready to introduce our treatments of aggregated modes, for which the aggregation matrix $\mathbf{P}$ may or may not be known.

Each aggregated mode is associated with two aspects and thus two factor matrices. {\em For mode subsampling:} two aspects would be subsampled independently based on the aforementioned categorical or continuous strategy. If they are connected with a known aggregation matrix $\mathbf{P}$, then $\mathbf{P}$ is subsampled accordingly along two dimensions; {\em For factor interpolation:} we follow the aforementioned strategies to interpolate two factors independently if the $\mathbf{P}$ is unknown; otherwise, we use the strategy to interpolate only the fine-granular factor and manually compute the other one. For example, at $d_t$ resolution, the aggregation matrix from $\mathbf{V}_{(d_{t})}$ to $\mathbf{Q}_{2(d_{t})}$, i.e., $\mathbf{P}_{2(d_{t})}$, is known. In this case, we would first interpolate a high resolution $\mathbf{V}_{(d_{t+1})}$ and then compute a high resolution $\mathbf{Q}_{2(d_{t+1})}$ based on $\mathbf{Q}_{2(d_{t+1})}=\mathbf{P}_{2(d_{t+1})}\mathbf{V}_{(d_{t+1})}$, instead of interpolating up from  $\mathbf{Q}_{2(d_{t})}$.

\section{Optimization Algorithms} \label{sec:optimization}
This section presents an effective alternating optimization algorithm, which is proposed for the coupled decomposition problem at each resolution. We first present the ALS algorithm and then formulate our quadratic subproblems.

\subsection{Coupled-ALS Algorithm}
\subsubsection{Alternating Least Squares (ALS)} To solve a low-rank CP approximation of a given tensor, ALS is one of the most popular algorithms \cite{sidiropoulos2017tensor}. It has been applied for both decomposition \cite{larsen2020practical,sidiropoulos2017tensor} and completion \cite{song2019tensor,10.1016/j.parco.2015.10.002}. However, our approach differs by considering a reconstructed (interim) tensor at each step in a  coupled decomposition setting. Below is the standard ALS algorithm.


Given a tensor $\mathcal{X}\in\mathbb{R}^{I_1\times I_2\times I_3}$, to find a rank-$R$ CP decomposition with factors $\mathbf{U}\in\mathbb{R}^{I_1\times R}$, $\mathbf{V}\in\mathbb{R}^{I_2\times R}$, $\mathbf{W}\in\mathbb{R}^{I_3\times R}$, the following loss function is minimized,
\begin{equation}
    \mathcal{L}(\mathcal{X};\mathbf{U,V,W}) = \|\mathcal{X}-\llbracket \mathbf{U,V,W}\rrbracket\|_F^2.
\end{equation}
ALS algorithm optimally solves for one of the factor matrices while keeping the others fixed and then cycles through all of the factor matrices. A sub-iteration of ALS for $\mathbf{U}$ is executed by solving  
\begin{equation}
    (\mathbf{V}\odot \mathbf{W})\mathbf{U}^\top =\mathbf{X}_1^\top,
\end{equation}
where $\mathbf{V}\odot \mathbf{W}\in\mathbb{R}^{I_2I_3\times R}$ is the Khatri-Rao product and $\mathbf{X}_1\in\mathbb{R}^{I_1\times I_2I_3}$ is the unfolded tensor along the first mode. 


\subsubsection{Coupled-ALS}
In our case, at each resolution $d_t$, we have a coupled tensor decomposition problem as defined in equation~\eqref{eq:dt_problem}. Similar to standard ALS, if we solve for one factor matrix while keeping all the other unknowns fixed, the problem reduces to a quadratic problem which can be solved optimally by coupling normal equations with respect to the factor at each granularity.

To reduce clutter, we elaborate on the full resolution form of Equation~\eqref{eq:original}. Specifically, we solve a least squares problem with the following equations.
\begin{itemize}
    \item For $\mathbf{U}$,  equations are derived using $\tilde{\mathbf{X}}_1$ and $\mathbf{C}^{(2)}_1$, which are the unfoldings of $\tilde{\mathcal{X}}$ and $\mathcal{C}^{(2)}$ along the first mode. Thus, we have the following two equations,
    \begin{align}
       (\mathbf{V}\odot \mathbf{W})\mathbf{U}^\top &=\tilde{\mathbf{X}}_1^{\top}, \label{eq:jacobi_example} \\
        \lambda_2(\mathbf{Q}_2\odot \mathbf{W})\mathbf{U}^\top &=\lambda_2\mathbf{C}_1^{(2)\top}. \label{eq:jacobi_second}
    \end{align}
    \item Similarly, for $\mathbf{V}$, equations are derived using $\tilde{\mathbf{X}}_2$ and $\mathbf{C}^{(1)}_2$.
    \item For $\mathbf{W}$, equations are derived using $\tilde{\mathbf{X}}_3$, $\mathbf{C}^{(1)}_3$ and $\mathbf{C}^{(2)}_3$.
    \item For $\mathbf{Q}_1$, equations are derived using $\mathbf{C}^{(1)}_1$.
    \item For $\mathbf{Q}_2$, equations are derived using $\mathbf{C}^{(2)}_2$ only. 
    However, $\mathbf{Q}_2$ and $\mathbf{V}$ satisfy the constraint $\mathbf{Q}_2=\mathbf{P_2V}$.
\end{itemize}


\subsection{Multi-resolution ALS}

\subsubsection{Alternating Optimization with Constraint.} Considering the constraint issue, e.g., between $\mathbf{V}$ and $\mathbf{Q}_2$, we use the equations to update one of them (the one corresponds to fine granular aspect) and later update the other one. For our example, $\mathbf{V}$ is updated by the linear equations, while $\mathbf{Q}_2$ is updated by the constraint $\mathbf{Q}_2=\mathbf{P}_2\mathbf{V}$ and new $\mathbf{V}$, instead of using linear equations\footnote{Note that, since these two factors are coupled, we could have transformed the $\mathbf{Q}_2$ equations to supplement $\mathbf{V}$ equations by pseudo inverse (i.e. inverse of $\mathbf{P}_2^T\mathbf{P}_2$), such that $\mathbf{V}$ would associate with three linear equations. However, the non-uniqueness of $\mathbf{P}_2^T\mathbf{P}_2$ makes the optimization unstable.}. The optimization flow follows the ALS routine, where during each iteration, we sequentially update $\mathbf{W,Q_1,U}$, $\mathbf{V}$ by the associated linear equations, and then manually update a new $\mathbf{Q_2}$. 

\subsubsection{Joint Normal Equations.} $\mathbf{W,Q_1,U}$, $\mathbf{V}$ are solved by using joint normal equations. For each of them, we consider stacking multiple linear systems and formulating joint normal equations with respect to each factor for the objective defined in equation~\eqref{eq:original}. Using $\mathbf{U}$ as the example, the normal equations are given as
\begin{equation}
    \begin{bmatrix}
    \mathbf{V}\odot \mathbf{W} \\
    (\mathbf{Q_2}\odot \mathbf{W})
    \end{bmatrix}^\top
    \begin{bmatrix}
    \mathbf{V}\odot \mathbf{W} \\
    \lambda_2(\mathbf{Q_2}\odot \mathbf{W})
    \end{bmatrix}\mathbf{U}^\top = \begin{bmatrix}
    \mathbf{V}\odot \mathbf{W} \\
    (\mathbf{Q_2}\odot \mathbf{W})
    \end{bmatrix}^\top
    \begin{bmatrix}
    \tilde{\mathbf{X}}_1^\top \\
    \lambda_2\mathbf{C}^{{(2)}\top}_1
    \end{bmatrix},
\end{equation}
where the components are weighted accordingly. The formulation could be further simplified as
\begin{align} \label{eq:joint}
        &((\mathbf{V}^\top\mathbf{V})*(\mathbf{W}^\top\mathbf{W}) + \lambda_2(\mathbf{Q}_2^\top\mathbf{Q}_2)*(\mathbf{W}^\top\mathbf{W}))\mathbf{U}^\top \notag\\
        &=(\mathbf{V}\odot \mathbf{W})^\top\tilde{\mathbf{X}}_1^{\top} + \lambda_2(\mathbf{Q}_2\odot \mathbf{W})^\top\mathbf{C}_1^{(2)\top}.
    \end{align}

This formulation is amenable to large scale problems as the right hand sides can exploit sparsity and the use of the matricized tensor times Khatri-Rao product (MTTKRP) kernel with each tensor. To avoid dense MTTKRP with $\tilde{\mathcal{X}}$, the implicit form given in Equation~\eqref{eq:implicit_form} can be used to efficiently calculate the MTTKRP by
\begin{align}
\label{eq:efficient_form}
        (\mathbf{V}\odot \mathbf{W})^\top\tilde{\mathbf{X}}_1^{\top}= (\mathbf{V}\odot \mathbf{W})^\top(\mathbf{M}_1 *& \mathbf{X}_1)^{\top} + (\mathbf{V}^\top \mathbf{V}^k * \mathbf{W}^\top \mathbf{W}^k)\mathbf{U}^{k\top}\notag\\ &-~(\mathbf{V}\odot \mathbf{W})^\top(\mathbf{M}_1 *\mathbf{R}_1)^{\top},
    \end{align}
where $\mathcal{R} = \llbracket \mathbf{U}^k,\mathbf{V}^k,\mathbf{W}^k \rrbracket$. The above simplification requires sparse MTTKRP and a sparse reconstruction tensor $\mathcal{M}*\mathcal{R}$, for which efficient programming abstractions are available~\cite{zhang2019enabling}.



\subsubsection{Solver Implementation} To solve the joint normal equations, we employ a two-stage approach.
\begin{itemize}
    \item {\bf Stage 1: } We start from the initialization given by the low resolution and apply the weighted Jacobi algorithm for $N$ (e.g., 5) iterations, which smooths out the error efficiently.
    \item {\bf Stage 2: } Starting from the $(N+1)-th$ iteration, we apply Cholesky decomposition to the left-hand side (e.g., Equation~\eqref{eq:joint}) and decompose the joint normal equations into solving two triangular linear systems, which could generate more stable solutions. Finally, the triangular systems are solved exactly.
\end{itemize}

In practice, we find that the combination of Stage 1 and Stage 2 works better than directly applying Stage 2 in light of both convergence speed and final converged results, especially when the number of observed entries is small. We perform ablation studies for this part. 

The {\bf time and space complexity} of \method is asymptotically equal to that of applying standard CP-ALS on  original tensor $\mathcal{X}$, which is $O(I_1I_2I_3R)$ and $O(I_1I_2I_3)$, respectively.

\section{Experiments}
The experiments provide a comprehensive study on \method. As an outline, the following evaluations are presented:
\begin{enumerate}
    \item ablation study on the different amount of partially observed data, compared to SOTA tensor completion and decomposition methods,
    \item performance comparison on tensor completion and downstream tensor prediction task, compared with baselines,
    \item ablation study on model initialization, compared with other initialization baselines.
\end{enumerate}

\begin{table*}[ht!] \small
\caption{Dataset Statistics}
\vspace{-2mm}
		\begin{tabular}{c|cccc|cccc} 
		\toprule 
		{\bf Name} & {\bf $1_{st}$-mode} & {\bf $2_{nd}$-mode} & {\bf $3_{rd}$-mode} & {\bf Sparsity} & {\bf Agg. $1_{st}$-mode} & {\bf Agg. $2_{nd}$-mode} & {\bf Agg. $3_{rd}$-mode} & {\bf Partial Obs.} \\
		\midrule
		HT & zipcode (1,200) & ICD-10 (1,000) & date (128) & 75.96\% & county (220) & CCS (189) & $--$ & 5\% \\
		GCSS & identifier (2,727) & keyword (422) & date (362) & 87.73\% & state (50) & $--$ & week (52) & 5\% \\
		\bottomrule 
		\multicolumn{6}{l}{
        \noindent~ * sparsity x\% means x\% of the entries are zeros in the original tensor.
        }
	    \end{tabular} \label{tb:dataset}
\end{table*}

\vspace{-2mm}
\subsection{Data Preparation}
\subsubsection{Two COVID-19 Related Databases}
Two COVID-related databases are considered in our evaluation: 
\begin{itemize}
    \item The {\em health tensor (HT)} is a proprietary real-world dataset, constructed from IQVIA medical claim database, including a {\em location by disease by time} tensor, counting the COVID related disease cases in the US from Mar. to Aug. 2020 during the pandemic.
    \item (ii) The {\em Google COVID-19 Symptoms Search (GCSS)}\footnote{https://pair-code.github.io/covid19\_symptom\_dataset/?country=GB} is a public dataset, which shows aggregated trends in Google searches for disease-related health symptoms, signs, and conditions, in format {\em location identifier by keyword by time}.
\end{itemize} 
\subsubsection{Data Processing}
For {\em HT} data, we use the most popular 1,200 zipcodes and 1,000 ICD-10 codes as the first two modes and use the first five months as the third mode (while the last month data is used for the prediction task), which creates $\mathcal{X}$. We aggregate the first mode by county and the second mode by CCS code separately to obtain two aggregated tensors, $\mathcal{C}^{(1)},\mathcal{C}^{(2)}$. We randomly sample 5\% elements from the tensor as the partial observation, $\mathcal{X}*\mathcal{M}$. For the {\em GCSS} data, we use all 2,727 location identifiers and 422 search keywords as the first two modes and a complete 2020 year span as the third time mode. To create two coarse viewed tensors, we aggregate on the state level in the first mode and week level for the third mode, separately. Assume 5\% of the elements are observed in {\em GCSS}. For both datasets, the first two modes are assumed categorical, while the last time mode is continuous.

In this experiment, we assume all the aggregation matrices are unknown. In the Appendix, we use synthetic data and {\em GCSS} to show that our model could achieve {\bf OracleCPD} level accuracy if some aggregation $\mathbf{P}$ is given. Basic data statistics are listed in Table~\ref{tb:dataset}.

\begin{table*}[h!] 
{\caption{Tensor Completion Results. {\normalfont Our \method outperforms the baselines significantly on PoF and CPU Time while having about the same space complexity. 
}}
\vspace{-1mm} \label{tb:recovery}
	\resizebox{0.95\textwidth}{!}{\begin{tabular}{c|ccc|ccc} \toprule 
	       \multirow{2}{*}{\bf Model} & \multicolumn{3}{c|}{\bf Health Tensor (HT)} & \multicolumn{3}{c}{\bf Google COVID-19 Symptoms Search (GCSS)} \\
	       \cmidrule{2-7}
	       & PoF & CPU Time & Peak Memory & PoF & CPU Time & Peak Memory\\
	       \midrule
	       BGD  & 0.4480 $\pm$ 0.0161 (7.155e-05) & 11,708.02 $\pm$ 98.39s & {\bf 10.27 GB} & 0.5552 $\pm$ 0.0098 (5.733e-06) & 26,074.53 $\pm$ 207.00s & 19.32 GB\\
	       B-PREMA  & 0.4678 $\pm$ 0.0176 (1.329e-04) & 11,614.75 $\pm$ 113.34s  & 11.04 GB  & 0.5923 $\pm$ 0.0045 (1.791e-06) &  26,201.12 $\pm$ 273.30s & 19.40 GB\\
	       CMTF-OPT  & 0.4815 $\pm$ 0.0245 (4.490e-04) & 10,159.50 $\pm$ 88.37s  & 10.53 GB & 0.5720 $\pm$ 0.0198 (9.394e-05) & 25,681.41 $\pm$ 141.30s & {\bf 19.25 GB}\\
	       \midrule
	       $\method_{both-}$  & 0.6109 $\pm$ 0.0234 (5.713e-02) & 1,437.20 $\pm$ 67.57s & {10.32 GB}  & 0.7508 $\pm$ 0.0063 (2.022e-01) & {\bf 1,821.84 $\pm$ 93.39s} & {19.75 GB}\\
	       $\method_{multi-}$  & 0.6432 $\pm$ 0.0144 (4.821e-01) & {\bf 1,375.32 $\pm$ 52.66s} & {10.31 GB} & 0.7560 $\pm$ 0.0059 (6.880e-01) & {1,874.18 $\pm$ 74.37s} & {20.09 GB}\\
	       $\method_{stage1-}$  & 0.6282 $\pm$ 0.0223 (1.880e-01) & 1,499.26 $\pm$ 59.66s & 10.55 GB & 0.7524 $\pm$ 0.0044 (2.274e-01) & 1,932.02 $\pm$ 25.25s & 20.33 GB\\
	       \method & {\bf 0.6520 $\pm$ 0.0113} & 1,501.43 $\pm$ 58.36s & 10.73 GB & {\bf 0.7579 $\pm$ 0.0049} & 1,929.49 $\pm$ 40.27s & 20.15 GB\\
			\bottomrule
			\multicolumn{7}{l}{
        \noindent~ * table format: mean $\pm$ standard deviation ($p$-value)
        }
\end{tabular}}}
\end{table*}

\vspace{-0.5mm}
\subsection{Experimental Setup}
\subsubsection{Baselines}
We include the following comparison models from different perspectives: (i) reference models: {\bf CPC-ALS} \cite{10.1016/j.parco.2015.10.002}: state of the art tensor completion model, only using partial observed data, {\bf OracleCPD} (CP decomposition of the original complete tensor with ALS); (ii) related baselines (they are gradient based): {\bf Block gradient descent (BGD)} \cite{xu2013block}, {\bf B-PREMA} \cite{almutairi2019prema}, {\bf CMTF-OPT} \cite{acar2011all}. They use both partial and coarse information; (iii) initialization models: {\bf MRTL} \cite{park2020multiresolution}, {\bf TendiB} \cite{almutairi2020tendi}, {\bf Higher-order SVD (HOSVD)}; (iv) our variants: {\bf $\method_{multi-}$} removes multiresolution module; {\bf $\method_{stage1-}$} removes stage 1 (Jacobi) in the solver; {\bf $\method_{both-}$} removes both multiresolution part and stage 1. More details are in the Appendix.

\vspace{-0.5mm}
\subsubsection{Tasks and Metrics} 
Evaluation metrics include: (i) {\em Percentage of Fitness (PoF)} \cite{song2019tensor}, suppose the target tensor is $\mathcal{X}$ and the low-rank reconstruction is $\mathcal{R}=\llbracket\mathbf{U},\mathbf{V},\mathbf{W}\rrbracket$, the relative standard error (RSE) and PoF are defined by
    \begin{equation}
        RSE = \frac{\|\mathcal{X}-\mathcal{R}\|_F}{\|\mathcal{X}\|_F}, \quad \quad PoF = 1-RSE;
    \end{equation}
(ii) {\em CPU Time}; (iii) {\em Peak Memory}.

We also consider {\bf (2) future tensor prediction}, which belongs to a broader spatio-temporal prediction domain \cite{wang2020deep}. We evaluate related to tensor-based baselines, while a thorough comparison with other models is beyond our scope. Specifically, we use Gaussian Process (GP) with radial basis function and white noise kernel along the date aspect to estimate a date factor matrix for the next month,
\begin{equation}
    \mathbf{W}_{future} = GP(\mathbf{W}).
\end{equation}
Then, the future disease tensor could be estimated by $\mathcal{R}_{future}=\llbracket\mathbf{U},\mathbf{V},\mathbf{W}_{future}\rrbracket$.
We use the next-month real disease tensor to evaluate the prediction results, and again, we use PoF as the metric.

\subsubsection{Hyperparameters}  For {\em HT} dataset, we use CP rank $R=50$, $20$ iterations at each low resolution, and $N=5$ in solver stage 1, within which $5$ Jacobi rounds are performaned per iteration. For {\em GCSS}, we use $R=20$, $10$ iterations at each low resolution, $N=5$, and $10$ Jacobi rounds per iteration. By default, we set $200$ iterations at the finest resolution for both datasets, which ensures the convergence. The parameter $\lambda_i$ is set to be $\lambda=e^{-\frac{i}{20}}$ at the $i$-th iteration, which varies the focus from coarse information to fine-granular information gradually. We publish our dataset and codes in a repository\footnote{https://github.com/ycq091044/MTC}.

\vspace{-0.5mm}
\subsection{Comparison with Reference Models}

\begin{figure}
    \centering
    \includegraphics[width=0.3\textwidth]{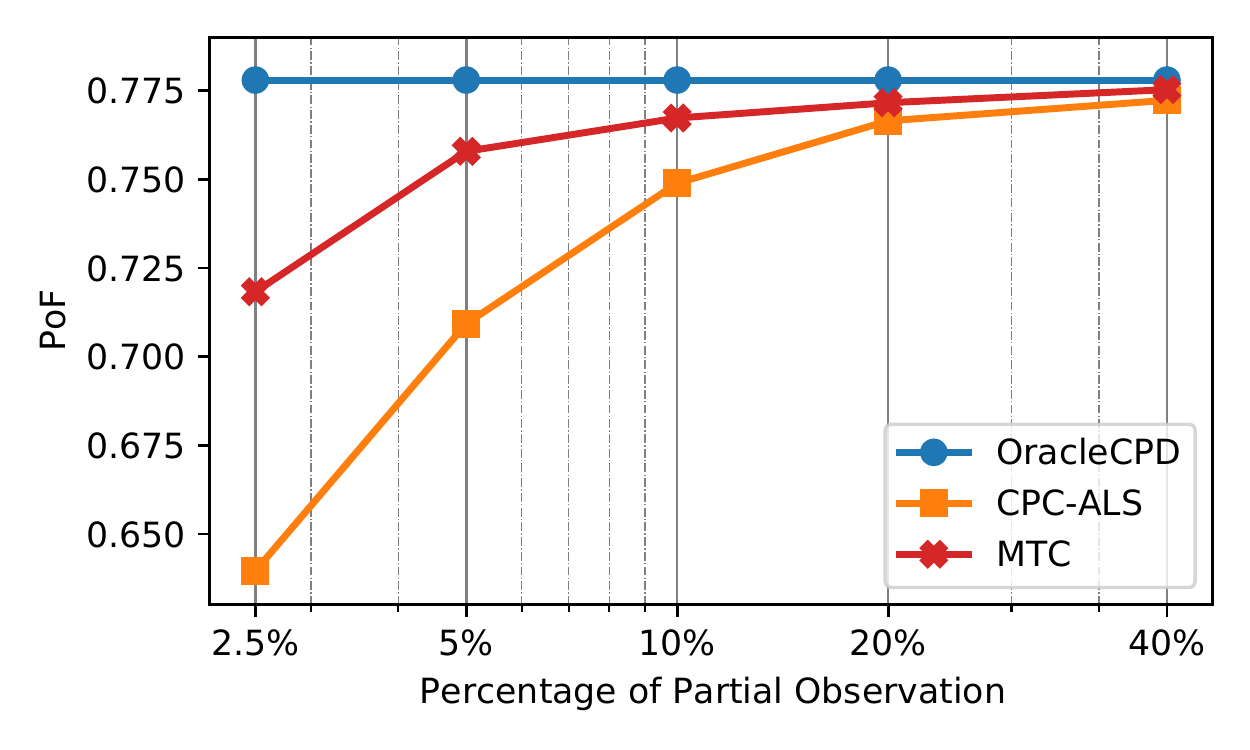}
    \vspace{-2mm}
    \caption{Ablation Study with Different Amount of Partially Observed Data on {\em GCSS}. {\normalfont Coarse information could improve tensor completion, especially when the partial data is not enough.}}
    \label{fig:reference}
\end{figure}

We first compare our \method with tensor decomposition model OracleCPD and a state of the art tensor completion model CPC-ALS on {\em GCSS}. The experiment shows the importance of different data sources and the model sensitivity w.r.t. the amount of partial data. Three models use the same CP rank $R=20$. The OracleCPD is implemented on the original tensor, which provides an ideal low-rank approximation. The CPC-ALS runs on the partially observed data only, while our \method uses both the partial and coarse tensors. 

We show the results in Figure~\ref{fig:reference}. When less partially observed data (e.g., 2.5\%) are presented, \method can achieve more (e.g., 7.5\%) improvement on PoF over CPC-ALS, with the help of coarse-granular data.  When we have more partial observations (e.g., 20\%), the gaps between three models become smaller (within 5\%).

\vspace{-1mm}
\subsection{Comparison on Tensor Completion}

We compare our \method and the variants with related baselines: BGD, B-PREMA, CMTF-OPT on the tensor completion task. The experiments are conducted 3 times with different random seeds. The results, standard deviations, and statistical $p$-values are shown in Table~\ref{tb:recovery}. We can conclude that the baselines are inferior to our proposed methods by around 15\% in fitness measure on both datasets.


The table also shows the memory requirements and executed time of each method.
\method reduces the time consumption significantly to about $\frac18$ of that of the baselines, since: (i) our iterative method incurs fewer floating-point arithmetics than computing the exact gradients from Equation~\eqref{eq:original}; (ii) gradient baselines require to perform expensive line search with multiple function evaluations.

\begin{figure}
    \centering
    \includegraphics[width=0.48\textwidth]{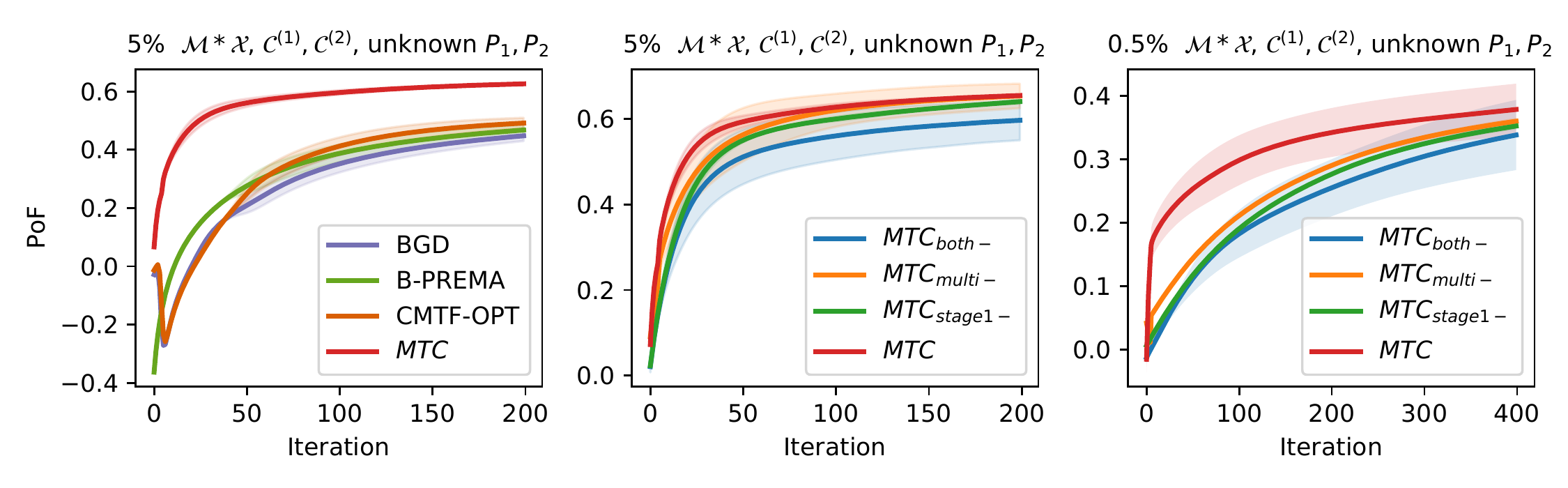}
    \vspace{-5mm}
    \caption{Optimization Landscape and Ablation Study on {\em HT}. {\normalfont Compared to gradient based baselines, \method provides a faster convergence rate and better converged results. With less partial observed information, \method becomes more advantageous.}}
    \vspace{-0.5mm}
    \label{fig:landscape}
\end{figure}

We also show the optimization landscape of all models and some ablation study results in Figure~\ref{fig:landscape}. We find that: (i) the multiresolution factorization and solver stage 1 are better used together; (ii) use of \method is more beneficial when having more missing values. 

\begin{figure*}
    \centering
    \includegraphics[width=0.95\textwidth]{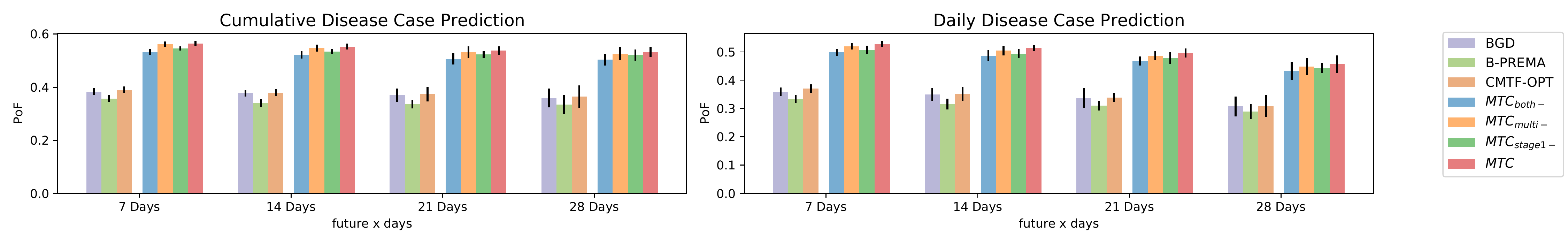}
    \vspace{-4mm}
    \caption{Predicting Future Disease Cases with the Learned Factors on {\em HT}. {\normalfont The results of tensor prediction accord with the tensor completion results, where our \method and the variants outperform the baselines by a great margin. In addition, we conjecture that there might be some underlying weekly patterns in the data, which explains why model performance decreases slowly with longer prediction windows.}}
    \label{fig:prediction_task}
\end{figure*}

\subsection{Comparison on Tensor Prediction}
After the optimization on {\em HT} dataset, we load the CP factors of each method and conduct the {\em tensor prediction task} as described in the experimental setup. Specifically, we perform cumulative and daily disease case prediction, where for one particular ICD-10 code at one particular zipcode location, we predict how many more disease cases will happen within the future $x$ days. 

The results are shown in Figure~\ref{fig:prediction_task}, where the dark bars are the standard deviations. The prediction results confirmed \method and its variants outperform the baselines by up to about 20\% on PoF.

\subsection{Ablation Study on Initializations}
To evaluate the multiresolution factorization part in \method, we compare with other tensor initialization methods. Specifically, we implement the initialization from (i) other baselines: MRTL, TendiB; (ii) a typical tensor initialization method: HOSVD; (iii) no initialization, $\method_{multi-}$. For a fair comparison, we use the same optimization method (in Section~\ref{sec:optimization}) after different initializations.

\begin{figure}[tbp!]
    \centering
    \includegraphics[width=0.48\textwidth]{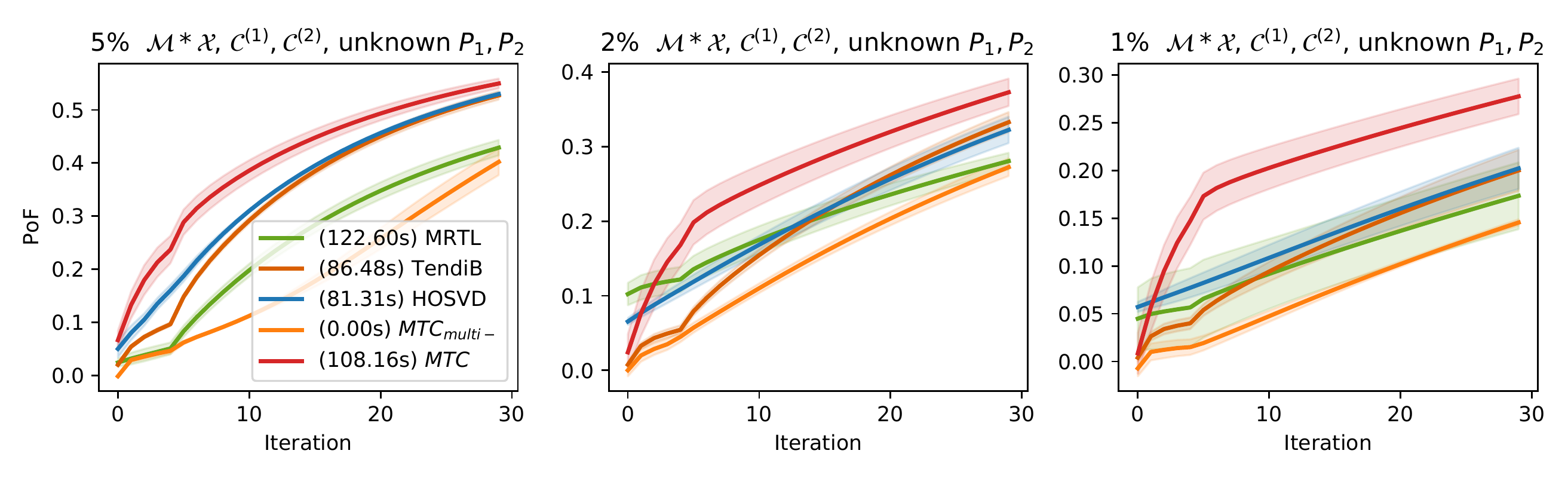}
    \vspace{-7mm}
    \caption{Initialization Comparison on {\em HT}. {\normalfont The multiresolution factorization in \method outperforms other initialization methods. With the multiresolution factorization module, the $5$ Jacobi iterations in stage 1 become very powerful, and then the stage 2 Cholesky solver continues to provide a stable and accurate estimation.}}
    \vspace{-2mm}
    \label{fig:initialization}
\end{figure}

Different amount of partial observations (i.e., different $\mathcal{M}$) are considered: 5\%, 2\% and 1\% on {\em HT} dataset. We report the first 30 iterations and the initialization time in Figure~\ref{fig:initialization}. Essentially, our \method outperforms all other models after 30 iterations. One interesting finding is that  TendiB and HOSVD start from the same PoF level regardless of the amount of partial observation since their initialization barely depends on the coarse views, i.e., $\mathcal{C}^{(1)}$ and $\mathcal{C}^{(2)}$. However, with only a few iterations, \method can surpass all baselines. 

\section{Conclusion}
This paper identifies and solves a new problem: {\em tensor completion from its partial and coarse observations.} We formulate the problem into a generic couple-ALS form and propose an efficient completion model, \method. The module combines multiresolution factorization and an effective optimization method as the treatment. We evaluate our \method on two COVID-19 databases and show noticeable performance gain over the baseline methods.

\section*{Acknowledgements}

This work was in part supported by the National Science Foundation award SCH-2014438, PPoSS 2028839, IIS-1838042, the National Institute of Health award NIH R01 1R01NS107291-01 and OSF Healthcare.

\bibliographystyle{ACM-Reference-Format}
\bibliography{sample-base.bib}

\appendix
\section{Basics of Tensor Computation}

\paragraph{Kronecker Product.} One important product for matrices is Kroncker product. For $\mathbf{A}\in\mathbb{R}^{I\times J}$ and $\mathbf{B}\in \mathbb{R}^{K\times L}$, their Kroncker product is defined by (each block is a scalar times matrix)
\begin{equation} \small
    \mathbf{A}\otimes\mathbf{B} = \begin{bmatrix}
    \mathbf{A}(1,1)\cdot \mathbf{B} & \mathbf{A}(1,2)\cdot \mathbf{B} & \cdots&\mathbf{A}(1,J)\cdot \mathbf{B} \\
    \mathbf{A}(2,1)\cdot \mathbf{B}&\mathbf{A}(2,2)\cdot \mathbf{B}&\cdots&\mathbf{A}(2,J)\cdot \mathbf{B}\\
    \vdots&\vdots&\cdots&\vdots\\
    \mathbf{A}(I,1)\cdot \mathbf{B}&\mathbf{A}(I,2)\cdot \mathbf{B}&\cdots&\mathbf{A}(I,J)\cdot \mathbf{B}
    \end{bmatrix}\in\mathbb{R}^{IK\times JL}. \notag
\end{equation}
\paragraph{Khatri–Rao Product.} 
Khatri-Rao product is another important product for matrices, specifically, for matrices with same number of columns. The Khatri-Rao product of $\mathbf{A}\in\mathbb{R}^{I\times J}$ and $\mathbf{B}\in \mathbb{R}^{K\times J}$ can be viewed as column-wise Kroncker product,
\begin{equation}
    \mathbf{A}\odot\mathbf{B} = \begin{bmatrix}
    \mathbf{A}(:, 1)\otimes\mathbf{B}(:, 1),~ \mathbf{A}(:, 2)\otimes\mathbf{B}(:, 2),~ \cdots,~ \mathbf{A}(:, L)\otimes\mathbf{B}(:, L)
    \end{bmatrix}, \notag
\end{equation}
where $\mathbf{A}\odot\mathbf{B}\in\mathbb{R}^{IK\times L}$.

\paragraph{Tensor Unfolding.} This operation is to matricize a tensor along one mode. For tensor $\mathcal{X}\in\mathbb{R}^{I_1\times I_2\times I_3}$, we could unfold it along the first mode into a matrix $\mathbf{X}_{1}\in\mathbb{R}^{I_1\times I_2I_3}$ (we use subscript notation). Specifically, each row of $\mathbf{X}_{1}$ is a vectorization of a slab in the original tensor; we have
\begin{equation}
    \mathbf{X}_1(i,j\times I_3+k) = \mathcal{X}(i,j,k). \notag
\end{equation}
Similarly, for the unfolding operation along the second or third mode, we have
\begin{align}
    \mathbf{X}_2(j,i\times I_3+k) &= \mathcal{X}(i,j,k)\in\mathbb{R}^{I_2\times I_1I_3}, \notag\\
    \mathbf{X}_3(k,i\times I_2+j) &= \mathcal{X}(i,j,k)\in\mathbb{R}^{I_3\times I_1I_2}. \notag
\end{align}

\paragraph{Hadamard Product.} The Hadamard product is the element-wise product for tensors of the same size. For example, the Hadamard product of two $3$-mode tensors $\mathcal{X},\mathcal{Y}\in\mathbb{R}^{I_1\times I_2\times I_3}$ is 
\begin{equation}
    \mathcal{Z} = \mathcal{X} * \mathcal{Y} \in\mathbb{R}^{I_1\times I_2\times I_3}. \notag
\end{equation}

\section{Solver Implementation}
We show the implementation of two stages in our solver.
\paragraph{Stage 1.}
Suppose that we want to solve for the factor matrix $\mathbf{U}$, given Equation~\eqref{eq:joint}. To simplify the derivation, we assume
\begin{align}
        \mathbf{A} &= (\mathbf{V}^\top\mathbf{V})*(\mathbf{W}^\top\mathbf{W}) + \lambda_2(\mathbf{Q}_2^\top\mathbf{Q}_2)*(\mathbf{W}^\top\mathbf{W}), \label{eq:A_def}\\
        \mathbf{B} &=(\mathbf{V}\odot \mathbf{W})^\top\tilde{\mathbf{X}}_1^{\top} + \lambda_2(\mathbf{Q}_2\odot \mathbf{W})^\top\mathbf{C}_1^{(2)\top}.  \label{eq:B_def}
    \end{align}

The Jacobi method first decomposes $\mathbf{A}$ into a diagonal matrix $\mathbf{D}$ and an off-diagonal matrix $\mathbf{R}$, such that $\mathbf{A}=\mathbf{D}+\mathbf{R}$. After moving the off-diagonal part to the right, {Equation~\eqref{eq:joint}} becomes
\begin{equation}
    (\mathbf{D}+\mathbf{R})\mathbf{U}^\top = \mathbf{B} \quad\Rightarrow\quad \mathbf{D}\mathbf{U}^\top = \mathbf{B} - \mathbf{R}\mathbf{U}^\top.
\end{equation}

It is cheap to take the inverse of the diagonal $\mathbf{D}$ (before the inverse, we add a small $\epsilon=$1e-5 to the diagonal values to improve numerical stability) on both sides. The Jacobi iteration is given by the following recursive equation,
\begin{equation} \label{eq:jacobi_recursive}
    \mathbf{U}^\top = \mathbf{D}^{-1}\mathbf{B} - \mathbf{D}^{-1}\mathbf{R}\mathbf{U}^\top ~ \Rightarrow ~ \mathbf{U}^{\top(k+1)} = \mathbf{D}^{-1}\mathbf{B} - \mathbf{D}^{-1}\mathbf{R}\mathbf{U}^{\top(k)}.
\end{equation}

The convergence property is ensured by the largest absolute value of eigenvalues (i.e., spectral radius) of $\mathbf{D}^{-1}\mathbf{R}$. If the the eigenvalues of $\mathbf{D}^{-1}\mathbf{R}$ are all in $(-1,1)$, then the iteration will converge. The iteration is enhanced by adding first-order momentum.

\begin{align} \label{eq:jacobi_example2}
    \mathbf{U}^{\top(k+1)} &= (1-w)\mathbf{U}^{\top(k)} + w(\underbrace{\mathbf{D}^{-1}\mathbf{B} - \mathbf{D}^{-1}\mathbf{R}\mathbf{U}^{\top(k)}}_{\mbox{the Jacobi term}}) \notag\\
    &= w\mathbf{D}^{-1}\mathbf{B} + (1-w\mathbf{D}^{-1}\mathbf{A})\mathbf{U}^{\top(k)}.
\end{align}
Now, the new spectral radius of {$1-w\mathbf{D}^{-1}\mathbf{A}_1$} is controlled by $w$.

\paragraph{Stage 2.} This stage uses Cholesky decomposition to solve the joint normal equation.
We use {\em scipy.linalg.solve\_triangular} to get the exact solution.

\section{Details in Experiments}
We list the GCSS data, processing files, the synthetic data, codes of our models, and all other baselines in an open repository.
	
\subsection{Baseline Implementation}
We provide more details about baselines used in the experiments:
\begin{itemize}
    \item {\bf CPC-ALS} is a SOTA tensor completion method. It could be integrated with the sparse structure on parallel machines. The implementation could refer to \cite{10.1016/j.parco.2015.10.002}. We provide our implementation in the link.
    \item {\bf OracleCPD} is implemented by standard CP decomposition on the original tensor, which provides the target low-rank factors. We implement it with our proposed solver. After convergence, the PoF results and objective functions are similar to the results provided by python Package {\em tensorly.decomposition.parafac}.
    \item {\bf Block Gradient Descent (BGD)} uses $\mathbf{Q}_1,\mathbf{Q}_2$ to replace $\mathbf{P_1U},\mathbf{P_2V}$ in Equation~\eqref{eq:original_problem}, and then apply the block gradient descent algorithm. The learning rate is chosen by a binary line search algorithm with depth 3.
    \item {\bf B-PREMA} \cite{almutairi2019prema}: this work is similar to {BGD} but with extra loss terms,
    \begin{equation}
        \beta_1\|\mathbf{1}^\top \mathbf{Q}_1 - \mathbf{1}^\top\mathbf{U}\|^2 + \beta_2\|\mathbf{1}^\top \mathbf{Q}_2 - \mathbf{1}^\top\mathbf{V}\|^2,
    \end{equation} ensuring that after the aggregation, the respective column sums of this two factors should be equal. Then, it also uses line search to find the learning rate.
    \item {\bf CMTF-OPT} \cite{acar2011all} is a popular model for coupled matrix and tensor factorization. We adopt it for coupled tensor factorization. Instead of alternatively optimizing factor matrices, it updates all factors at once. We use gradient descent to implement it, and the learning rate is also selected by line search.
    \item {\bf MRTL} \cite{park2020multiresolution} designs a multiresolution tensor algorithm to build good initialization for accelerating their classification task. We adopt their tensor initialization method and compare it to our proposed multiresolution factorization. 
    \item {\bf TendiB} \cite{almutairi2020tendi} provides a heuristic initialization trick for the coupled factorization task. We only adopt its initialization based on our problem setting. In our case, we implement it as follows. We apply CPD on $\mathcal{C}^{(1)}$ to estimate $\mathbf{Q}_1$,  $\mathbf{V}$ and $\mathbf{W}$ first, and then we freeze $\mathbf{W}$ and use it to estimate $\mathbf{U}$ and $\mathbf{Q}_2$ by applying another CPD on $\mathcal{C}^{(2)}$. We have also tried to interpolate $\mathbf{U}$ and $\mathbf{Q}_2$ from the estimated $\mathbf{Q}_1$ and $\mathbf{V}$. However, it performs poorly. Consequently, we adopte the first implementation.
    \item {\bf Higher-order singular value decomposition (HOSVD)} \cite{kolda2009tensor} is commonly used for initializing tensor CP-ALS decomposition. It first alternatively estimates orthogonal factors by higher-order orthogonal iteration (HOOI) and then applies ALS on the core tensor. The final initialization is the compositions of the orthogonal factor and core ALS factor. We adopt the HOSVD method as another initialization baseline.
\end{itemize}

\subsection{Tasks and Metrics} 
We provide more details and rationality on the metrics.
\begin{itemize}
    \item {\bf Percentage of Fitness (PoF).} Several performance evaluation metrics have been introduced in a survey paper \cite{song2019tensor}, and PoF is the most common one. For this metric, we also perform a two-tailed Student T-test and showe the $p$-value in Table~3.
    \item {\bf CPU Time}. We make the methods execute for the same number of iterations and then compare their running time. The running time includes data loading, initialization (in our methods, it also includes the expense on low resolutions) and excludes the metric computation time.
    \item {\bf Peak Memory}. To evaluate space complexity, we record the memory usage during the optimization process and report the peak memory load.
\end{itemize}

\subsection{Implementation and Hyperparameters} All the experiments are implemented by {\em Python 3.7.8, scipy 1.5.2 and numpy 1.19.1} on Intel Cascade Lake Linux platform with 64 GB memory and 80 vCPU cores. By default, we conduct each experiment three times with different random seeds.

Since the CPD accepts an scaling invariance property, we apply the following trick after each iteration. For example, a set of factors $\{\mathbf{U,V,W}\}$ would be identical to another factor set $\{2\mathbf{U}, 0.5\mathbf{V},\mathbf{W}\}$ in terms of both the completion task or the prediction task. During the implementation of all models, we have adopted the following rescaling strategy to minimize the float-point round-off error: 
\begin{align}
    \mathbf{U}(:, i) &\leftarrow f_i\cdot \frac{\mathbf{U}(:, i)}{\|\mathbf{U}(:, i)\|}, \\
    \mathbf{V}(:, i) &\leftarrow f_i\cdot \frac{\mathbf{V}(:, i)}{\|\mathbf{V}(:, i)\|}, \\
    \mathbf{W}(:, i) &\leftarrow f_i\cdot \frac{\mathbf{W}(:, i)}{\|\mathbf{W}(:, i)\|},
\end{align}
where 
\begin{equation}
    f_i = \left(\|\mathbf{U}(:, i)\|\|\mathbf{V}(:, i)\|\|\mathbf{W}(:, i)\|\right)^{\frac13},~\forall i=[1..R],
\end{equation}
which essentially equilibrates the norms of the factors of each
component (we also do it for $\mathbf{Q}_2$ and $\mathbf{Q}_3$).

\subsection{Experiments on Synthetic Data} The experiments on synthetic data verify our model in the ideal low-rank case.
\paragraph{Synthetic Tensor.} The data is generated by three uniformly randomed rank-$10$ factors with $125$ as the mode size consistently, $\mathbf{U,V,W}\in\mathbb{R}^{125\times 10}$. We sort each column of $\mathbf{V,W}$, to make sure the mode smoothness. Thus, a tensor $\mathcal{X}=\llbracket\mathbf{U,V,W}\rrbracket\in\mathbb{R}^{125\times 125 \times 125}$ is constructed, assuming that the first mode is categorical and the second/third modes are continuous. We also generate two aggregation matrix $\mathbf{P}_1,\mathbf{P}_2\in\{0,1\}^{12\times 125}$ for the first and second modes, separately, to obtain coarse tensors, $\mathcal{C}^1=\mathcal{X}\times_1\mathbf{P}_1$, $\mathcal{C}^2=\mathcal{X}\times_2\mathbf{P}_2$.

\paragraph{Verification on Synthetic Data.}
On the synthetic data, \method is compared with three variants. We find that with only $2\%\sim3\%$ partially observed elements, all models can achieve near $100\%$ PoF, which is an exact completion of the low-rank synthetic $\mathcal{X}$. We thus further reduce the partial observation rate and consider three scenarios: (i) only $1\%$ observed data; (ii) $1\%$ observation with $\mathcal{C}^1,\mathcal{C}^2$, where $\mathbf{P}_1$ is unknown and $\mathbf{P}_2$ is known; (iii) $1\%$ observation with $\mathcal{C}^1,\mathcal{C}^2$, where $\mathbf{P}_1,\mathbf{P}_2$ are both unknown. 

\begin{figure}[hbp!]
    \centering
    \includegraphics[width=0.46\textwidth]{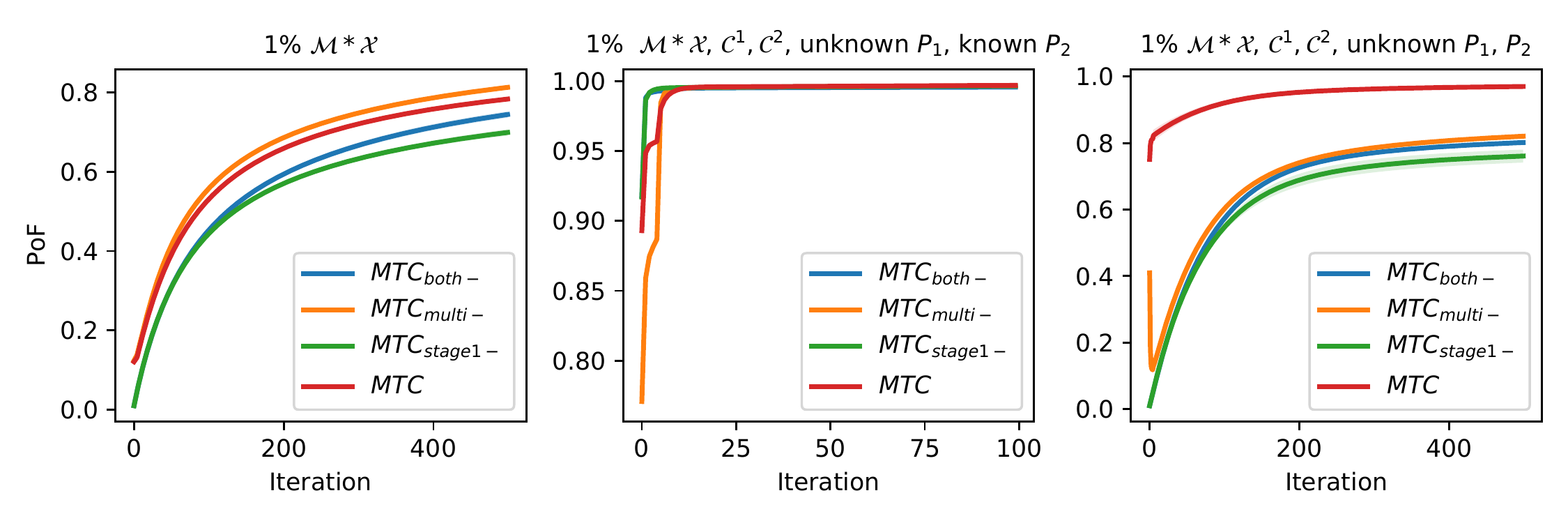}
    \vspace{-5mm}
    \caption{Comparison on Synthetic Data}
    \vspace{-4mm}
    \label{fig:synthetic}
\end{figure}

Results are shown in Figure~\ref{fig:synthetic}. Basically, we conclude that: (i) our model is more advantageous with coarse level information; (ii) the multiresolution factorization and stage 1 are better used together. 

\subsection{Additional Experiments on {\em GCSS} with Known Aggregation}
The main paper experiments do not use any aggregation information since, in those scenarios, the differences between model variants are clearer. This appendix supplies one aggregation matrix along the third model and shows that all our variants can achieve a similar tensor completion performance as OracleCPD.

Specifically, on {\em GCSS} dataset, we keep other settings unchanged and supply the aggregation matrix (from date to week) along the third mode. All variants converge quickly, and thus we only run the experiments for 50 iterations. The results are shown below.

\begin{table}[ht!] \small
\caption{Results with Known Aggregation on {\bf GCSS}}
\vspace{-1mm}
		\begin{tabular}{c|cccc|cccc} 
		\toprule 
		{\bf Model} & {\bf PoF}\\
		\bottomrule 
		{OracleCPD} & 0.7780 $\pm$ 0.0012  \\
		$\method_{both-}$ & 0.7706 $\pm$ 0.0058 \\
		$\method_{multi-}$ & 0.7694 $\pm$ 0.0025 \\
		$\method_{stage1-}$ & 0.7704 $\pm$ 0.0037 \\
		$\method$ & 0.7708 $\pm$ 0.0026 \\
		\bottomrule
	    \end{tabular}
\end{table}




\end{document}